\documentclass[reqno, 11pt]{amsart} 
\usepackage{amsfonts, amsmath, amssymb, amsthm}
\usepackage{mathrsfs}

\usepackage{amsfonts,amsmath,latexsym,verbatim,amscd,mathrsfs,color,array}
\usepackage[colorlinks=true]{hyperref}
\usepackage{amsmath,amssymb,amsthm,amsfonts,graphicx,color}
\usepackage{amssymb}
\usepackage{pdfsync}
\usepackage{epstopdf}
\usepackage{cite}
\usepackage{graphicx}
\usepackage{amsmath}
\usepackage{amsfonts}%
\providecommand{\U}[1]{\protect\rule{.1in}{.1in}}

\numberwithin{equation}{section}

\newtheorem{theorem}{Theorem}[section]

\newtheorem{proposition}[theorem]{Proposition}
\newtheorem{corollary}[theorem]{Corollary}
\newtheorem{remark}[theorem]{Remark}

\newtheorem{lemma}[theorem]{Lemma}

\newtheorem{definition}[theorem]{Definition}

\newtheorem*{nonumtheorem}{Theorem A}

\newcommand{\smalltext}[1]{{\fontsize{10}{12}\selectfont #1}}

\def\AR{{\mathcal A}} 
\def\BR{{\mathcal B}}

\def\DR{{\mathcal D}}
\def\ER{{\mathcal E}}
\def\FR{{\mathcal F}}

\def\JR{{\mathcal J}}

\def\LR{{\mathcal L}}

\def\NR{{\mathcal N}}

\def\UR{{\mathcal U}}

\def\WR{{\mathcal W}}

\def\ZR{{\mathcal Z}}

\def\LB{\mathbb{L}}

\def\rd{\mathrm d}   
\def\R{\mathbb{R}^n}

\def\hk{\hat K}
\def\hw{\hat W}
\def\hz{\hat Z}
\def\he{\hat E}
\def\hu{\hat U}

\def\vt{\vec{t}}
\def\vn{\vec{n}}
\def\ti{\tilde{I}}
\def\bq{\boldsymbol{q}}
\def\bp{\boldsymbol{p}}
\def\bbf{\boldsymbol{f}}
\def\bg{\boldsymbol{g}}
\def\bla{\boldsymbol{\lambda}}
\def\hbp{\hat{\boldsymbol{p}}}
\def\bc{\boldsymbol{c}}
\def\hbc{\hat{\boldsymbol{c}}}

\def\T{\text{T}}

\def\tw{\tilde{W}}

\def\omei{\Omega^\mathcal{I}}
\def\omem{\Omega^\mathcal{M}}
\def\omeo{\Omega^\mathcal{O}}

\begin{document}

	\title[Prescribed scalar curvature problem]{Infinitely many  solutions for the prescribed scalar curvature problem with volcano-like curvature}
	
	\author{Tuoxin Li}
	\address[Tuoxin Li]{Department of Mathematics, Chinese University of Hong Kong, Shatin, NT, Hong Kong}
	\email{txli@math.cuhk.edu.hk}
	
	\author{Juncheng Wei}
	\address[Juncheng Wei]{Department of Mathematics, Chinese University of Hong Kong, Shatin, NT, Hong Kong}
	\email{wei@math.cuhk.edu.hk}

	\author{Haidong Yang}
	\address[Haidong Yang]{School of Mathematics and Statistics, Central South
		University, Changsha, Hunan, PR China}
	\email{hdyang@csu.edu.cn}

	\maketitle

	\begin{abstract}
		In this paper, we consider the following prescribed scalar curvature problem:
		\begin{equation*}
			-\Delta u = K(x) u^{\frac{n+2}{n-2}}, \quad u>0\quad\hbox{in}\quad \mathbb{R}^n, \quad
			u \in D^{1,2}(\mathbb{R}^n),
		\end{equation*}
		where $K(x)$ is a volcano-like positive function such that
		$$ K(x)= K(r_0)- c_0 | |x|- r_0|^m + O( | |x|- r_0|^{m+\theta}),\quad r_0- \delta <|x| <r_0+\delta$$
		with $K(r_0), c_0, \delta>0, \theta >2, \min \{\frac{n-2}{2}, 2\} < m< n-2$.
		We first prove the existence of infinitely many positive solutions. A consequence of our proof yields that  the infinitely many solutions constructed in \cite{WY} are non-degenerate in the whole $D^{1, 2}(\mathbb{R}^{n})$ space. To our knowledge, it seems to be the first result of infinitely many solutions of prescribed scalar curvature problem when the potential function $K(x)$ is not radial. Our non-degeneracy results are also more complete and improve the result in \cite{GMPS}. 
	\end{abstract}

	\section{Introduction}\label{intro}
	~~~~The prescribed scalar curvature problem is a fundamental topic in differential geometry and geometric analysis. It asks for $n$-sphere $(\mathbb{S}^{n}, g_0)$ and fixed smooth function $\tilde{K}$, if there exists a conformally invariant metric $g_0$ such that the scalar curvature becomes $\tilde{K}$. Mathematically, this problem can be formulated as seeking solutions of the following equation:
	\begin{equation}\label{psp}
		\left\{
		\begin{array}{l}
			-\Delta_{\mathbb{S}^{n}} u + \frac{n(n-2)}{4} u = \tilde{K} u^{\frac{n+2}{n-2}} \quad \text{on} \; \mathbb{S}^{n}, \\
			u > 0.
		\end{array}
		\right.
	\end{equation}
	
	Using the stereo-graphic projection, the prescribed scalar
	curvature problem \eqref{psp} can be can be reduced to the following equation:
	\begin{equation}\label{eq u ori}
		\left\{
		\begin{array}{l}
			-\Delta u = K(x) u^{\frac{n+2}{n-2}}, u>0 \text{ in } \mathbb{R}^n, \\
			u \in D^{1,2}(\mathbb{R}^n).
		\end{array}
		\right.
	\end{equation}
	In the last three decades, enormous efforts have been devoted to the study of \eqref{eq u ori}. Many existence and non-existence results have been obtained  under various conditions of $K(x)$. In the simplest case, when $K(x)=K(|x|)\equiv 1$, \eqref{eq u ori} turns to
	\begin{equation}\label{eq yamabe}
		\Delta u+u^{\frac{n+2}{n-2}}=0,	\ u>0 \text{ in } \R,
	\end{equation}
	which is also called Yamabe equation. It is well-known that all the solutions to \eqref{eq yamabe} has the form (which are called Aubin-Talenti bubbles)
	\begin{equation}\label{AT bubbles}
		U_{Q,\Lambda}(x)=c_n\frac{\Lambda^{\frac{n-2}{2}}}{(1+\Lambda^2|x-Q|^2)^\frac{n-2}{2}},\quad c_n=({n(n-2)})^{\frac{n-2}{4}}
	\end{equation}
	for any $\Lambda >0$ and $Q \in\mathbb{R}^n$. We refer to \cite{Talen,CGS, GNN}. Thanks to above celebrated results to Yamabe equation, researchers have developed a variety of solutions for $K(x)$ under specific conditions. As $K(x)$ is a perturbation of constant, Cao, Noussair and Yan in \cite{CNY} proved that there are two or many solutions to \eqref{eq u ori}. Li constructed infinitely many multi-bump solutions of \eqref{eq u ori} in \cite{Li1,Li2,Li3} when $K(x)$ is a positive and periodic function and $n\geq 3$. In cases where $K(x)$ is 1-periodic in its first $l$ variables and satisfies others nature conditions near its critical point, Li, the second author and Xu place the centers of bumps in some lattices in $\mathbb{R}^l$ and proved the existence of this special multi-bump solutions in \cite{LWX}. In addition to the ones mentioned above, we refer to \cite{A1,Y1,Y2,Y3,CY,CFM,CM1,CI,CGYang93,ChenLin97,CL1,CL2,DLY,DingNi,GGZ,J1,LPZ,M1,Q1,XuWei,Yan} for the progress of the prescribed curvature problem in the fields of geometry and partial differential equations.
	
	In this paper, we are interested in the multi-bump solutions whose energy can be made arbitrary large. The second author and Yan in \cite{WY} proved the existence of the above type of solutions as $K(x)$ is radially symmetric. More specifically, assume that $K(x)$ satisfies the following condition.
	
	\medskip
	$(\mathbf{K}_0):$  $K(x)=K(r)$ is radially symmetric, and there is $r_0>0$ such that
	\begin{equation}\label{K0}
		K(r) = K(r_0) - c_0|r - r_0|^m + O(|r - r_0|^{m+\theta}), \quad r \in (r_0 - \delta, r_0 + \delta),
	\end{equation}
	where $K(r_0)>0, c_0 > 0 $, $ \theta > 0 $ are some constants, and $ m $ satisfies $ m \in [2, n-2) $.
	\medskip
	
	Authors in \cite{WY} proved the following theorem:
	\begin{nonumtheorem}\label{thm wei yan}
		Suppose that $n\ge 5$. If $K(r)$ satisfies $(\mathbf{K}_0)$, then problem \eqref{eq u ori} has infinitely many
		non-radial solutions.
	\end{nonumtheorem}
	
	By taking $k$, the number of the bubbles, as the parameter and applying finite reduction method, authors of \cite{WY} constructed a sequence of $k$-bump solution $\{u_k\}$. It is worth emphasizing that the approach of taking the number of bubbles as a free parameter, which was first introduced in \cite{WY}, has been extensively adopted in the construction of infinitely many solutions for various partial differential equations—particularly in the context of elliptic problems lacking small parameters. We refer to \cite{DMFP,MPW, SWY, WY1} for application of this method and \cite{BT,GPY,KW,LLY,LNW,NT} for the study of multi-bump solutions of other nonlinear equations.

	On the other hand, the proof presented in \cite{WY} relies in an essential way on the radial symmetry of the function $K(x)$. In fact, with this symmetry, one can build a $k$-bump solution for an arbitrary large $k$ with a $k$-dyadic symmetry, and the problem is then reduced to adjusting two parameters: the height of each bubble and the radius of the $k$-polygon. There is no doubt that the radial symmetry property of the potential function plays an important role in studying the existence of solutions to elliptic equations. The analysis becomes significantly more challenging when the potential function lacks symmetry, we refer to \cite{AWY,AW,CE,DWY} for the study of Schr\"odinger equations with non-symmetric potentials. In contrast, the analogous issue for the prescribed curvature problem remains considerably less explored. It is then natural to ask following question.
	
	\medskip
	
	\noindent
	\textbf{(Q): Whether Theorem A remains true when the symmetry requirement of $K(x)$ is removed? }
	
	\medskip
	
	Note that we cannot restrict our configuration to any symmetry class, as the function $K(x)$ lacks radial symmetry. Consequently, we are faced with a large number of bumps, and thus a large set of parameters that require careful adjustment. This introduces substantial difficulties in the construction compared to \cite{WY}. In this paper, we give an affirmative answer to the above question. To proceed, we first impose the following assumptions on $K(x)$.
	
	\medskip
	
	$(\mathbf{K}_1):$  There is $r_0>0$ such that
	\begin{equation}\label{K1}
		K(x)=K(r_0)-c_0||x|-r_0|^m+O(||x|-r_0|^{m+\theta}),\  |x|\in (r_0-\delta, r_0+\delta),
	\end{equation}
	where $K(r_0)>0, c_0 > 0 $, $ \theta > 2 $ are some constants, and $ m $ satisfies 
	$\min\{\frac{n-2}{2},2\}<m<n-2$.
	
	\medskip
	
	To state our result, we introduce the following function:
	\begin{equation*}
		g(x)=\sum_{j=1}^{\infty}\frac{1-\cos (jx)}{j^{n}}.
	\end{equation*}
	Up to a normalization constant, when $n$ is even, $g(x)$ is related to the Fourier series of the Bernoulli polynomial, and when $n$ is odd $g(x)$ is related to the Fourier series of the Euler polynomial. We refer to \cite{AB} for further details.
	
	The first result of this paper is stated as follows.
	\begin{theorem}\label{thm main}
		Let $n\ge 5$ and $K$ satisfies $(\mathbf{K}_1)$. Assume
		\begin{equation}\label{condi g}
			g''(x)<\frac{n-2}{n-1}\frac{(g'(x))^2}{g(x)}, \text{ for all } x\in (0,\pi).
		\end{equation}
		Then \eqref{eq u ori} admits infinitely many non-radial positive solutions, whose energy can be made arbitrarily large.
	\end{theorem}

	\begin{remark}
		Condition \eqref{condi g} was first proposed in \cite{MW} for proving the non-degeneracy of sign-changing solutions to Yamabe problem, constructed by del Pino, Musso, Pistoia and Pacard \cite{DMFP}. Meanwhile, it can be learned from \cite{MW} that \eqref{condi g} should be true for $n\ge 4$ since Musso and the second author had checked \eqref{condi g} numerically, up to dimension $n\le 48$. We refer to \cite{MW}  for more details.
	\end{remark}
	
	\begin{remark}
		This type of multi-bump solutions have also been constructed for the nonlinear Schr\"{o}dinger equation. We refer the reader to \cite{DWY} for the standard case and \cite{AWY} for the fractional case. We would like to make some more comments on the condition $(\mathbf{K}_1)$. As we will see in later computation,  the condition on $m$ ensures that the interaction of bubbles  is small, and hence the approximate solutions is good enough; while we set the condition on $\theta$ to make the remainder term small enough. From our proof it is not difficult to see that \eqref{K1} in $(\mathbf{K}_1)$ can be improved to be
		\begin{equation*}
			K(x)=K(r_0)-c_0||x|-r_0|^m+\tilde{c}||x|-r_0|^{m+\tilde{\theta}}+O(||x|-r_0|^{m+\theta}),\  |x|\in (r_0-\delta, r_0+\delta),
		\end{equation*}	
		where $\tilde{c}\in\R$ and $\tilde{\theta}>0$.
	\end{remark}
	
	A second objective of this article is to establish the non-degeneracy of the multi-bump solutions to equation \eqref{eq u ori}. When applying the reduction method, it is essential to verify that the building blocks are non-degenerate. This condition ensures their structural rigidity under the invariant transformations of the equation, a fundamental property in the analysis of partial differential equations. We refer the readers to \cite{KM1,KM2,KKK} for more applications of non-degeneracy property in the study of elliptic and parabolic equations.
	
	In a previous work, Guo, Musso, Peng and Yan in \cite{GMPS} showed that if $K(x)$ satisfies $(\mathbf{K}_0)$ with $m=2$ and $\theta=1$, and in addition
	\begin{equation}\label{condi GMPS}
		\Delta K-(\Delta K-\frac{1}{2}(\Delta K)')r\neq 0 \text{ at } r=r_0,
	\end{equation} 
	then the $k$-bump solution constructed in \cite{WY} are non-degenerate  the symmetric space $H_s$ (This concept of non-degeneracy here will be made precise in Proposition \ref{prop nondegen sym 1}). Here $H_s$ is defined by
	\begin{equation}\label{def Hs}
		\begin{aligned}
			H_{s}:=  \bigg\{ u : & u \text{ is even in } x_h, h = 2, \cdots, n,\\
			& u(r \cos \theta, r \sin \theta, y'') = u\left(r \cos\left(\theta + \frac{2\pi j}{k}\right), r \sin\left(\theta + \frac{2\pi j}{k}\right), y''\right)
			\bigg\}.
		\end{aligned}
	\end{equation} 
	
	Their proof is mainly based on some local Pohozaev identities. It is well-known that the non-degeneracy and 	the local uniqueness are two closely related properties of solutions, so similar method is often used to prove the local uniqueness of solutions. We refer to \cite{GLW, GJNT, GPY, LY} and the references therein.
	
	In this paper, we aim to obtain the non-degeneracy of multi-bump solutions constructed in Theorem A in a general class of functions. We first introduce our definition of non-degeneracy. Let $v$ be a solution to \eqref{eq u ori}. The linearized operator around $v(x)$ is defined by
	\begin{equation*}
		L_v =\Delta+\frac{n+2}{n-2}Kv^{\frac{4}{n-2}},
	\end{equation*}
	and the null space of $L_v$ is defined by
	\begin{equation}\label{def Zv}
		\mathcal{Z}_v :=\{f\in D^{1, 2} (\mathbb{R}^{n}): L_v f=0\}.
	\end{equation}
	We observe that if $K$ is radially symmetric, then equation \eqref{eq u ori} is invariant under rotation in $\R$. This invariance is reflected in the kernels of $L_v$, that is, we have $\bar{\mathcal{Z}}_{v}\subset \mathcal{Z}_{v}$, where
	\begin{equation}\label{def bar Zv}
		\bar{\mathcal{Z}}_{v} = \operatorname{span} \left\{
		\left(x_{j} \partial_{x_{i}}-x_{i} \partial_{x_{j}}\right)v, \ 1 \leq j < i \leq n
		\right\}.
	\end{equation}
	Inspired by\cite{DKM}, we now introduce our definition of non-degeneracy.
	\begin{definition}\label{def nondegen}
		A solution $v$ of \eqref{eq u ori} is said to be non-degenerate if $\mathcal{Z}_v=\bar{\mathcal{Z}}_{v}$.  
	\end{definition}
	To state our result, we impose the following condition on $K(x)$.
	
	\medskip
	$(\mathbf{K}_2):$  $K(x)=K(r)$ is radially symmetric and bounded, and there is $r_0>0$ such that
	\begin{equation*}
		K(r) = K(r_0) - c_0|r - r_0|^m + O(|r - r_0|^{m+\theta}), \quad r \in (r_0 - \delta, r_0 + \delta),
	\end{equation*}
	where $K(r_0)>0, c_0 > 0 $, $ \theta > 2 $ are some constants, and $ m $ satisfies $\min\{\frac{n-2}{2},2\}<m<n-2$.
	\medskip
	
	Then our non-degeneracy result can be stated as follows:
	\begin{theorem}\label{thm nondegen}
		Let $n\ge 5$ and $K(x)$ satisfy $(\mathbf{K}_2)$. Assume \eqref{condi g} holds. Then there exist infinitely many multi-bump solutions of \eqref{eq u ori} which are non-degenerate in the sense of Definition \ref{def nondegen}.
	\end{theorem}
	Using the same technique, we can also establish the non-degeneracy result in $H_s$:
	\begin{corollary}\label{thm nondegen sym}
		Let $n\ge 5$ and $K(x)$ satisfies $(\mathbf{K}_0)$. Then there exist infinitely many multi-bump solutions of \eqref{eq u ori} which are non-degenerate in $H_s$.
	\end{corollary}
	We remark that the condition here is the same as in Theorem A, and we do not need any other additional assumptions like \eqref{condi GMPS}. This improves Theorem 1.1 of Guo, Musso, Peng and Yan in \cite{GMPS}.
	
	Before closing this section, we shall summarize the main difficulties and novelties developed in this study.\par
	(\textbf{1}) To establish Theorem \ref{thm main}, we employ an intermediate Lyapunov–Schmidt reduction method to address equation  \eqref{eq u ori}. A similar procedure has been implemented in the context of Schr\"odinger equations by del Pino, the second author, and Yao in \cite{DWY}. Nevertheless, our setting presents two major challenges. First, the algebraic decay of the Aubin-Talenti bubbles and their interactions with all other bubbles introduce significant complexity into the reduced problem. Second, in addition to the parameters governing bubble locations, an extra parameter is required to characterize the height of each bubble. These parameters possess distinct scaling properties (see \eqref{def q norm}), necessitating more refined error estimates and a careful spectral analysis of the matrix $T$ arising in the reduced problem (see \eqref{def T}). This constitutes the main contribution of the present work. We believe that the analytical approach developed here can be adapted to construct other types of multi-bump solutions in which the number of bumps tends to infinity.\par
	(\textbf{2})  As we mentioned earlier, the non-degeneracy of solutions has many applications in the study of partial differential equations. Here we shall mention the recent work \cite{SWY} of Sun, the second author and Yang in solving the Brezis' open question. In the proof process \cite{SWY}, a very crucial step is that they chose a new type of solutions, which was called 'crown solutions', with sign-changing properties and fast decay rate as the building blocks. This enabled the finite-dimensional reduction method to be applied and eventually solved the open problem proposed by Brezis and Nirenberg \cite{BN} in 1983.
	Moreover, the methods developed in this paper are expected to be applicable to studying the non-degeneracy of multi-bump solutions in other nonlinear equations. In fact, one can also employ the same procedure to establish the non-degeneracy of the multi-bump solutions constructed in \cite{WY1}, as both the solutions and the potential function exhibit a parallel structure and invariance properties.
	\par
	(\textbf{3})
	The results of this work are unified in their demonstration that the constructive method of \cite{WY} remains effective for problems with "almost radial" potentials. It is important to note that the heightened computational complexity inherent in this generalization gives rise to the following open question:
	
	\medskip
	\noindent\textbf{Open Question:  Whether there exist non-degenerate multi-bump solutions to \eqref{eq u ori} as $K(x)$ is non-radial? }
	\medskip
	
	This paper is organized as follows. In Section \ref{prelim}, we briefly describe the construction of multi-bump solutions and present useful estimates. We perform the first Lyapunov-Schmidt reduction in Section \ref{reduction} and calculate projections along the kernels in Section \ref{proj}. The final proof of Theorem \ref{thm main} is presented in Section \ref{pf of thm main}. We place the relevant content of non-degeneracy results in Section \ref{Pf non-degen}. The Appendix is devoted to the detailed proofs of several computations.

	\section{Preliminaries}\label{prelim}
	~~~~Inspired by \cite{WY}, when considering $k$-bump solutions, we use the scaling $u(x)\rightarrow \mu^{-\frac{n-2}{2}}u(\frac{x}{\mu})$ to transform \eqref{eq u ori} to 
	\begin{equation}\label{eq u}
		\left\{
		\begin{array}{l}
			-\Delta u = \hat K(x) u^{\frac{n+2}{n-2}}, \quad u > 0, \; x \in \mathbb{R}^n, \\
			u \in D^{1,2}(\mathbb{R}^n),
		\end{array}
		\right.  
	\end{equation}
	where the scaling parameter $\mu:=k^{\frac{n-2}{n-2-m}}$ and $\hk(x):=K(\frac{x}{\mu})$. In the rest of paper, we denote $p=\frac{n+2}{n-2}$, and we always consider the above rescaled equation instead of the original one \eqref{eq u ori}. In addition, we also assume $K(r_0)=1$ without loss of generality when $K$ satisfies $(\mathbf{K}_0)$, $(\mathbf{K}_1)$ or $(\mathbf{K}_2)$.
	
	We next provide a brief description of the solutions to be constructed in Theorem \ref{thm main}, along with an outline of the principal ideas underlying the construction procedure. The solutions are obtained as small perturbations of a sum of Aubin–Talenti bubbles defined in \eqref{AT bubbles}, centered at a set of carefully selected points on $\mathbb{R}^2\times \{\mathbf{0}\}$.
	
	Let the locations of bubbles be given by $Q_j$, $j=1,\dots,k$. Then we define
	\begin{equation*}
		U_{j}(x):=U_{Q_j,\Lambda_j}(x)=c_n\frac{\Lambda_j^{\frac{n-2}{2}}}{(1+\Lambda_j^2|x-Q_j|^2)^\frac{n-2}{2}},\quad W=\sum_{j=1}^k U_j.
	\end{equation*}
	Here $\Lambda_j$ has the form $\Lambda_j=\Lambda+\lambda_j$, where $\Lambda$ is a fixed constant to be determined later, and $\lambda_j$'s are small perturbations.
	
	To describe further the configuration space of $Q_j$'s, we define
	\begin{equation*}
		Q_{j,\alpha}=(R\cos\theta_{j,\alpha}, R\sin\theta_{j,\alpha},\mathbf{0}),
	\end{equation*}
	where $\theta_{j,\alpha}=\frac{2\pi(j-1)}{k}+\alpha$, and $R$ is some constant to be determined later. As we mentioned before, here $\alpha$ is the parameter that is used to deal with the degeneracy due to rotations. If  $\alpha=0$, then the $Q_{j,\alpha}$'s are the points used in \cite{WY}. It is obvious that the parameter $\alpha$ plays no role in the construction as $K(x)$ is radially symmetric. However, it is very important in our construction as we will see in Section \ref{pf of thm main}. 
	
	For the constant $R$, it has the form $R=\mu r_0+R_0$. Let $\theta=\frac{2\pi}{k}$. Then $R_0$ and $\Lambda$ should satisfy the following balancing condition:
	\begin{align}
		&c_0\mu^{-m}(\gamma_2\Lambda^{2-m}R_0+\frac{\gamma_4\Lambda^{-m}}{R})+\gamma_1\Lambda^{2-n}(2R)^{1-n}\sum_{l=1}^{k-1} \sin^{2-n}(l\theta)=0,\label{R0 La 1}	  \\
		&\gamma_1\Lambda^{2-n}(2R)^{2-n}\sum_{l=1}^{k-1}\sin^{2-n}(l\theta)-2c_0\mu^{-m}\gamma_3\Lambda^{-m}=0,\label{R0 La 2}	
	\end{align}
	where $\gamma_1, \gamma_2, \gamma_3$ and $\gamma_4$ are given as the following:
	\begin{equation}\label{ga 1}
		\gamma_1=(n-2)c_n^{p+1}\int_{\R}\frac{1}{(1+|y|^2)^\frac{n+2}{2}}\rd y.
	\end{equation}
	\begin{equation}\label{ga 2}
		\gamma_2=(n-2)mc_n^{p+1}\int_{\R}
		\frac{|y_1|^m}{(1+|y|^2)^{n+1}}\rd y.
	\end{equation}
	\begin{equation}\label{ga 3}
		\gamma_3=\frac{m}{p+1}\int_{\R} |y_1|^m
		U^{p+1}(y)\rd y.
	\end{equation}
	\begin{equation}\label{ga 4}
		\gamma_4=\frac{n-2}{2}mc_n^{p+1}\int_{\R}
		\frac{|y_1|^m|y|^2-|y_1|^{m+2}}{(1+|y|^2)^{n+1}}\rd y.
	\end{equation}
	Substituting \eqref{R0 La 2} into \eqref{R0 La 1}, we get
	\begin{equation*}		\gamma_2\Lambda^{2}RR_0+\gamma_4+\gamma_3=0.
	\end{equation*} 
	Then it is easy to check that it has a unique solution in the regime $|R_0|\le \frac{C}{\mu}$. This balancing condition can be understood as a condition such that the approximation $W$ is very close to a genuine solution. In other words, this is also the condition that the bubbles should satisfy (at the main order) when $K(x)$ is radially symmetric (see \eqref{r La0} and \eqref{r La1}).
	
	Since $K(x)$ can be viewed as a small perturbation of a radial function, we also choose our $Q_j$, the location of the bubble, to be small perturbation of $Q_{j,\alpha}$. To be more specific, for $f_j, g_j\in \R$, we define
	\begin{equation}\label{def Qj}
		Q_j=Q_{j,\alpha}+f_j\vn_j+g_j\vt_j,
	\end{equation}
	where
	\begin{equation*}
		\vn_j=(\cos\theta_{j,\alpha}, \sin\theta_{j,\alpha}, \mathbf{0}),\quad \vt_j=(-\sin\theta_{j,\alpha}, \cos\theta_{j,\alpha},\mathbf{0}).
	\end{equation*}
	We remark that $f_j$ and $g_j$ measure the displacement around $Q_{j,\alpha}$ in the normal and tangential directions respectively.

	To quantitatively measure the bubbles $U_j$ and their interactions, we introduce a vector-valued parameter  $\bq=(\bla, \bbf, \bg)\in \mathbb{R}^{k}\times\mathbb{R}^{k}\times \mathbb{R}^{k}$. For $K$ satisfies  $(\mathbf{K}_1)$, we fix a small positive number $\tau$ such that 
	\begin{equation}\label{def tau}
		\tau<\min\{\theta-2, \frac{m-2}{2}, \frac{2m+2-n}{m}\}.
	\end{equation}
	Then we fix a small constant $0<\tau_1<\frac{1}{2}\tau$, and define our parameter space as
	\begin{equation}\label{def Xi}
		\Xi:=\{\bq=(\bla, \bbf, \bg)\in \mathbb{R}^{3k}:\|\bq\|_\Xi\le 1\},
	\end{equation}
	where
	\begin{equation}\label{def q norm}
		\|\bq\|_\Xi:=\left\{ \begin{aligned}&\mu (d^{\tau_1}\|\bla\|_{\infty}+d^{\frac{3}{2}\tau_1}\|\bbf\|_{\infty})+d^{\tau_1}\left(\|\bg\|_{\infty}+\sup_{j\neq l}\Big|\frac{g_j-g_l}{\sin(j-l)\theta}\Big|\right), &n\ge 6,\\
			&\mu (\|\bla\|_{\infty}+\|\bbf\|_{\infty})+d^{\tau_1}\|\bg\|_{\infty}, &n=5.
		\end{aligned}\right.
	\end{equation}
	
	We prove Theorem \ref{thm main} by showing the following result.
	\begin{theorem}\label{thm u}
		Let $n\ge 5$ and $K$ satisfies $(\mathbf{K}_1)$. Assume \eqref{condi g} holds. Then for sufficiently large $k$, \eqref{eq u} has at least two solutions of the form
		\begin{equation*}
			u(x)=W(x)+\phi(x),
		\end{equation*}
		where $\alpha\in [0,2\pi)$ and the parameter $\bq\in \Xi$. Moreover, $\|\phi\|_{L^\infty} \to 0$ as $k\to \infty$.
	\end{theorem}

	We next present some preliminary facts and some useful estimates. Let $d_{jl}=|Q_{j,\alpha}-Q_{l,\alpha}|$ and $\rho_{jl}=|Q_j-Q_{l}|$. Then for $\bq\in\Xi$, it is easy to see that $\rho_{jl}=d_{jl}+O(\frac{1}{\mu})$. In particular, for $d:=|Q_{j,\alpha}-Q_{j+1,\alpha}|=2R\sin\theta$ and $\rho:=\min\limits_{j}|Q_j-Q_{j+1}|$, we have
	\begin{equation}
		\rho=d+O(\frac{1}{\mu}).
	\end{equation}
	From this we see that $\rho$ is very close to $d$, so in many cases we don't need to distinguish them and we prefer to use $d$ in this paper for simplicity.
	
	As we will see in the next section, when we apply the Lyapunov-Schmidt reduction, technical difficulties arise as the number of bumps goes to infinity. To deal with these difficulties, we divide $\R$ into $3k$ parts. For $j=1, \dots, k$, let 
	\begin{equation}\label{def omega j}
		\Omega_j:=\{ x\in\R: |x-Q_{j,\alpha}|=\min_{1\le l\le k}\{|x-Q_{j,\alpha}|\}\}.
	\end{equation}
	Then we define
	\begin{equation}\label{def omei}
		\omei_j:=\Omega_j\cap\{x\in\R:  |x-Q_{j,\alpha}|\le \frac{d}{2}\},
	\end{equation}
	\begin{equation}\label{def omem}
		\omem_j:=\Omega_j\cap\{x\in\R:  |x-Q_{j,\alpha}|\in (\frac{d}{2},\frac{1}{2}\delta\mu)\},
	\end{equation}
	\begin{equation}\label{def omeo}
		\omeo_j:=\Omega_j\cap\{x\in\R:  |x-Q_{j,\alpha}|\ge \frac{1}{2}\delta\mu\}.
	\end{equation}
	Then $\R=\bigcup_{j=1}^k (\omei_j\cup \omem_j\cup \omeo_j)$ and the interior of these sets are disjoint. Sometimes we don't need to distinguish $\omem_j$ and $\omeo_j$, so we denote $\Omega_j'=\omem_j\cup\omeo_j$ for  $j=1, \dots, k$.
	
	We remark that though it seems more natural to use the points $Q_{j}$ to define above sets, it is simpler to use $Q_{j,\alpha}$ to define them, and their difference is indeed small enough. For example, since $|Q_{j,\alpha}-Q_{j}|\le C$, we have
	\begin{equation}\label{diffenrence}
		\big|\omei_j\setminus B_{\frac{d}2}(Q_j)\big|+ \big|B_{\frac{d}2}(Q_j)\setminus\omei_j \big|\le Cd^{n-1}.
	\end{equation}
	
	Now we give two estimates that will be frequently used in this paper. We fix a constant $\sigma>1$ which is close enough to  $1$. Then for $x\in \Omega_i$ and $s>2$, we have
	\begin{equation*}
		\begin{aligned}
			\sum_{j=1}^k\frac{1}{1+|x-Q_j|^{s}}\le& \frac{1}{1+|x-Q_i|^{s}}+ \frac{C}{1+|x-Q_i|^{s-\sigma}}\sum_{j\neq i}\frac{1}{|Q_i-Q_j|^{\sigma}}\\
			\le& \frac{1}{1+|x-Q_i|^{s}}
			+\frac{Cd^{-\sigma}}{1+|x-Q_i|^{s-\sigma}}.
		\end{aligned}
	\end{equation*}
	In particular, we have
	\begin{equation}\label{est W i}
		\sum_{j=1}^k\frac{1}{1+|x-Q_j|^{s}}\le \left\{\begin{aligned}
			&\frac{C}{1+|x-Q_i|^{s}}, &x\in \omei_{i},\\
			&\frac{Cd^{-\sigma}}{1+|x-Q_i|^{s-\sigma}}, &x\in \Omega_{i}'.
		\end{aligned}\right.
	\end{equation}
	On the other hand, we fix a number $\sigma_1\in (\frac{n-2-m}{n-2},1)$ which is close enough to $\frac{n-2-m}{n-2}$. Then since
	\begin{equation*}
		\sum_{j\neq i}\frac{1}{|Q_i-Q_j|^{\sigma_1}}\le Cd^{-\sigma_1}\sum_{j=1}^k\frac{1}{j^{\sigma_1}}\le C,
	\end{equation*}
	we also have
	\begin{equation}\label{est lose}
		\sum_{j=1}^k\frac{1}{1+|x-Q_j|^{s}}\le
		\frac{C}{1+|x-Q_i|^{s-\sigma_1}},\ x\in \Omega_{i}'.
	\end{equation}

	We close this section by providing some estimates of convolution type with polynomial decays, which will be used to analyze the interactions between bumps. The proofs can be found in \cite{B1}.
	
	\begin{equation}\label{E1}
		\int_{\R} pU_{j}^{p-1}\nabla U_j U_{l}\rd x
		=\gamma_1 \Lambda_j^{-\frac{n-2}{2}}\Lambda_l^{-\frac{n-2}{2}}\frac{Q_j-Q_l}{|Q_j-Q_j|^n}+O(\rho_{jl}^{-n-1}).
	\end{equation}	
	
	\begin{equation}\label{E2}
		\begin{aligned}
			\int_{\R} pU_{j}^{p-1}Z_{j,0} U_{l}\rd x
			=&-\frac{1}{2}\gamma_1
			\Lambda_j^{-\frac{n}{2}}\Lambda_l^{-\frac{n-2}{2}}|Q_j-Q_l|^{2-n}
			+O(\rho_{jl}^{-n}).
		\end{aligned}
	\end{equation}
	
	\begin{equation}\label{E3}
		\begin{aligned}
			\int_{\R} pU_{j}^{p-1}\nabla U_{j}Z_{l,0} \rd x=\frac{n-2}{2}\gamma_1\Lambda_j^{-\frac{n-2}{2}}\Lambda_l^{-\frac{n}{2}}\frac{Q_j-Q_l}{|Q_j-Q_l|^n}+O(\rho_{jl}^{-n-1}).
		\end{aligned}
	\end{equation}
	
	\begin{equation}\label{E4}
		\begin{aligned}
			\int_{\R} p(p-1)U_{j}^{p-2}\nabla U_{j}Z_{j,0}U_l \rd x=\frac{n-2}{2}\gamma_1\Lambda_j^{-\frac{n}{2}}\Lambda_l^{-\frac{n-2}{2}}\frac{Q_j-Q_l}{|Q_j-Q_l|^n}+O(\rho_{jl}^{-n-1}).
		\end{aligned}
	\end{equation}
	
	\begin{equation}\label{E5}
		\begin{aligned}
			\int_{\R} pU_{j}^{p-1}Z_{j,0} Z_{l,0}\rd x
			=\frac{n-2}{4}\gamma_1
			\Lambda_j^{-\frac{n}{2}}\Lambda_l^{-\frac{n}{2}}|Q_j-Q_l|^{2-n}
			+O(\rho_{jl}^{-n}).
		\end{aligned}
	\end{equation}

	\section{The Lyapunov-Schmidt Reduction}\label{reduction}
	~~~~The aim of this section is to complete a first reduction procedure. To do this, we first introduce some definitions of weighted norms employed in this paper.
	
	Inspired by \cite{DSWduke}, to measure errors and the interactions among various bubbles more accurately, we define following weighted functions. Let $\sigma>1$ be a fixed constant which is close enough to $1$. For $Q_j$ defined in  \eqref{def Qj}, we define
	\begin{equation}\label{def of V}
		\mathcal{V}(x)=\sum_{i=1}^{k}\frac{d^{-(n-2)}}{1+|x-Q_i|^{4}}\chi_{\{x\in \omei_i\}}+\sum_{i=1}^{k}
		\frac{d^{-\max\{p,2\}\sigma}}{1+|x-Q_i|^{n+2-\max\{p,2\}\sigma}}\chi_{\{x\in \Omega'\}},
	\end{equation}
	and
	\begin{equation}\label{def of W}
		\mathcal{W}(x)=\sum_{i=1}^{k}\frac{d^{-(n-2)}}{1+|x-Q_i|^{2}}\chi_{\{x\in \omei_i\}}+\sum_{i=1}^{k}
		\frac{d^{-\max\{p,2\}\sigma}}{1+|x-Q_i|^{n-\max\{p,2\}\sigma}}\chi_{\{x\in \Omega'\}}.
	\end{equation}
	Then we define our norms as
	\begin{equation}\label{def of n1}
		\|h\|_{*}:=\sup \mathcal{V}(x)^{-1}|h(x)|,
	\end{equation}
	and
	\begin{equation}\label{def of n2}
		\|h\|_{**}:=\sup \mathcal{W}(x)^{-1}|h(x)|.
	\end{equation}
	
	For $u=W+\phi$, solving \eqref{eq u} is equivalent to solving
	\begin{equation}\label{eq phi}
		L(\phi)+E+N(\phi)=0,
	\end{equation}
	where
	\begin{equation}\label{def L}
		L(\phi)=\Delta\phi+p\hk W^{p-1}\phi, 
	\end{equation}
	\begin{equation}\label{def E}
		E=\Delta W+\hk W^p=(\hk-1)W^p+(W^p-\sum_{j=1}^{k}U_{j}^p),
	\end{equation}
	and
	\begin{equation}\label{def N}
		N(\phi)=\hk\left((W+\phi)^p-W^p-pW^{p-1}\phi\right).
	\end{equation}
	
	Instead of solving \eqref{eq phi} directly, we will solve a projected version of it. To be more precise, we define
	\begin{equation}\label{def Zj}
		Z_{j,0}=\frac{\partial U_{j}}{\partial\Lambda_j},\ 
		Z_{j,1}=\nabla U_{j}\cdot \vn_j,\ 
		Z_{j,2}=\nabla U_{j}\cdot \vt_j,
	\end{equation}
	and $\tilde{W}=\frac{1}{p}U_j^p$.
	Then direct computation shows that
	\begin{equation*}
		\frac{\partial\tw}{\partial\bq}=(U_1^{p-1}Z_{1,0},\dots, U_k^{p-1}Z_{k,0}, -U_1^{p-1}Z_{1,1},\dots, -U_k^{p-1}Z_{k,1}, -U_1^{p-1}Z_{1,2},\dots, -U_k^{p-1}Z_{k,2}).
	\end{equation*}
	
	Now we can state our main result in this section.
	\begin{proposition}\label{prop nonline proj}
		Suppose that $n\ge5$ and ${K}$ satisfies $(\mathbf{K}_1)$. Let $Z_{j,l}$ be given in \eqref{def Zj}. Then there is a positive integer $k_0$ such that for all $k\ge k_0$, every $\alpha\in \mathbb{R}$ and $\bq\in\Xi$, there exist a unique function $\phi\in W^{2,2}(\mathbb{R}^n)\cap \mathcal{B}_{*}$ 
		where
		\begin{equation}\label{def B*}
			\mathcal{B}_{*}:=\{\phi\in L^{\infty}(\mathbb{R}^n):\ \|\phi\|_{*}\le C_0\}
		\end{equation}
		for some large constant $C_0$, and a unique multiplier $\bc=(c_{j,l})\in\mathbb{R}^{3k}$ such that
		\begin{equation}\label{eq phi proj}
			\left\{\begin{aligned}
				&L(\phi)+E+N(\phi)=\bc\cdot\frac{\partial\tw}{\partial\bq},\\
				&\int_{\R}\phi U_{j}^{p-1}Z_{j,l}\rd x=0, j=1,\dots,k, l=0,1,2.
			\end{aligned}\right.
		\end{equation}
		Moreover, the map $(\alpha, \bq)\rightarrow \phi(x;\alpha, \bq)$ is of class $C^1$ with
		\begin{equation}\label{est of dependence}
			R^{-1}\|\frac{\partial \phi}{\partial \alpha}\|_{*}
			+\|\frac{\partial \phi}{\partial \bq}\|_{*}
			\le C.
		\end{equation}		
	\end{proposition}
	The proof of Proposition \ref{prop nonline proj} is now standard. We will first study the invertibility of the linear operator, and then apply the fixed point theorem to the nonlinear equation.
	
	\subsection{Linear Theory}
	Our main result in this subsection can be formulated as the following.
	\begin{proposition}\label{prop line theo}
		Under the assumption of Proposition \ref{prop nonline proj}. There is a positive integer $k_0$ such that for all $k\ge k_0$, every $\alpha\in \mathbb{R}$,  $\bq\in\Xi$, and all $h$ with $\|h\|_{**}<\infty$, there exist a unique function $\phi\in W^{2,2}(\R)$ and a unique multiplier $\bc\in\mathbb{R}^{3k}$ such that
		\begin{equation}\label{line theo}
			\left\{\begin{aligned}
				&L(\phi)=h+\bc\cdot\frac{\partial \tilde{W}}{\partial\bq}\\
				&\int_{\R}\phi U_{j}^{p-1}Z_{j,l}\rd x=0, j=1,\dots,k, l=0,1,2.
			\end{aligned}\right.
		\end{equation}
		Moreover, we have estimates
		\begin{equation}\label{priori est of phi and c}
			\|\phi\|_*\le C\|h\|_{**},\ \|\bc\|_{\infty}\le d^{-(n-2)}\|h\|_{**}.
		\end{equation}	
	\end{proposition}
	
	The estimates of $L(Z_{j,l})$ play an important role in the proof of Proposition \ref{prop line theo}, so we first establish the following lemma:
	\begin{lemma}\label{coro norm L}
		For sufficiently large $k$ and $l=0,1,2$, there exists a constant $C$ such that
		\begin{equation}\label{norm sum L Zjl}
			\Big\|\sum_{j=1}^k L(Z_{j,l})\Big\|_{**}\le C.
		\end{equation}
	\end{lemma}
	\begin{proof}
		By definition,
		\begin{equation*}
			L(Z_{j,l})=\Delta Z_{j,l}+p\hk W^{p-1}Z_{j,l}
			=(\hat{K}-1) W^{p-1}Z_{j,l}+
			(W^{p-1}-U_{j}^{p-1})Z_{j,l}.
		\end{equation*}
		In the following, we focus on the case $l=1,2$, since the estimate for $L(Z_{j,0})$ is similar.
		
		If $x\in \Omega_i'$ for some $i$, by \eqref{est W i}, we have
		\begin{equation}\label{est L Zjl i'}
			|L(Z_{j,l})|\le\frac{Cd^{-(p-1)\sigma}}{|x-Q_i|^{4-(p-1)\sigma}}\frac{1}{|x-Q_j|^{n-1}}.
		\end{equation}
		If $x\in \omei_i$, then similar to \eqref{est W i}, we have
		\begin{equation*}
			|(W^{p-1}-U_{j}^{p-1})Z_{j,l}|\le
			\left\{\begin{aligned}
				&\frac{C}{1+|x-Q_i|^4} 	\frac{1}{1+|x-Q_j|^{n-1}}, &i\neq j,\\
				&\frac{C}{d^{n-2}}\frac{1}{1+|x-Q_j|^{5}}, &i=j,
			\end{aligned}\right.
		\end{equation*}
		and
		\begin{equation*}
			|(\hat{K}-1) W^{p-1}Z_{j,l}|\le
			\left\{\begin{aligned}
				&\frac{C\mu^{-m}|x-Q_i|^m}{1+|x-Q_i|^4} 	\frac{1}{|x-Q_j|^{n-1}},&i\neq j,\\
				&C\mu^{-m}\frac{|x-Q_j|^m}{1+|x-Q_j|^{n+3}}, &i=j.
			\end{aligned}\right.
		\end{equation*}
		Then since $m<n-2$ and $\mu^m\sim d^{n-2}$, we obtain
		\begin{equation}\label{est L Zjl 0}
			|L(Z_{j,l})|\le \left\{\begin{aligned}
				&\frac{C}{1+|x-Q_i|^4} 	\frac{1}{1+|x-Q_j|^{n-1}}, &x\in \omei_i, i\neq j,\\
				&\frac{C}{d^{n-2}}\frac{1}{1+|x-Q_j|^{5}}, &x\in \omei_j.
			\end{aligned}\right.
		\end{equation}
		
		For $Z_{j,0}$, the same argument as in \eqref{est L Zjl i'} and \eqref{est L Zjl 0} shows that
		\begin{equation}\label{est L Zj0}
			|L(Z_{j,0})|\le \left\{\begin{aligned}
				&\frac{Cd^{-(p-1)\sigma}}{|x-Q_i|^{4-(p-1)\sigma}}\frac{1}{|x-Q_j|^{n-2}}, &x\in \Omega_i',\\
				&\frac{C}{1+|x-Q_i|^4} 	\frac{1}{1+|x-Q_j|^{n-2}}, &x\in \omei_i, i\neq j,\\
				&\frac{C}{d^{n-2}}\frac{1}{1+|x-Q_j|^{4}}, &x\in \omei_j.
			\end{aligned}\right.
		\end{equation}
		
		Taking summation in $j$, we obtain \eqref{norm sum L Zjl} from the definition of $\|\cdot \|_{**}$.
	\end{proof}

	\begin{proof}[Proof of Proposition \ref{prop line theo}]
		It is enough to prove a priori estimates, since the rest of the argument is the same as the proof of Proposition 4.1 in \cite{DMcv}. To be more specific, we will establish the following:
		
		\noindent \textbf{Claim:} Assume $\phi=\phi_k$ solves \eqref{line theo} for $h=h_k$. If $\|h_k\|_{**}\to 0$ as $k\to\infty$, then so does $\|\phi_k\|_*$.
		
		We prove this claim by contradiction. If the claim is not true, then we can find a sequence $\phi_{k}$ solves \eqref{line theo} for $h_{k}$ with $\|\phi_k\|_{*}=1$ and $\|h_k\|_{**} \to 0$. For simplicity, we drop the subscript $k$ in the rest of the argument. 
		
		We first establish a priori estimates of $\bc=\{c_{j,l}\}$. For $j=1,\cdots, k$, $l=0,1,2$, multiplying the first equation of \eqref{line theo} by  $Z_{j,l}$ and integrating over $\mathbb{R}^n$, we have
		\begin{equation}\label{multiply against kernel}
			\int_{\mathbb{R}^n}L(\phi)Z_{j,l}\rd x=\int_{\mathbb{R}^n}hZ_{j,l}\rd x
			+\sum_{s=0}^2\sum_{t=1}^k c_{t,s}\int_{\mathbb{R}^n}U_{t}^{p-1}Z_{t,s}Z_{j,l}\rd x.
		\end{equation}
		For the coefficients of $c_{t,s}$ on the right-hand side, standard arguments gives that
		\begin{equation}\label{almost ortho}
			\int_{\R}U_{t}^{p-1}Z_{j,l}Z_{t,s}\rd x=\left\{
			\begin{aligned}
				&	\int_{\R}U_{t}^{p-1}Z_{j,l}^2\rd x, t=j, s=l,\\
				&0, t=j, s\neq l,\\
				&O(d_{tj}^{2-n}), t\neq j, sl=0,\\
				&O(d_{sj}^{1-n}), t\neq j, sl\neq 0.
			\end{aligned}
			\right.
		\end{equation}
		Therefore the matrix $M_0$ defined by $(M_0)_{tj,sl}=\int_{\R}U_{t}^{p-1}Z_{t,s}Z_{j,l}\rd x$ is strictly diagonally dominant.
		
		We next estimate the other two terms in \eqref{multiply against kernel}. 
		For $n\geq 6$,	using  \eqref{est W i}, \eqref{est L Zjl 0} and the fact that $\|\phi\|_*=1$, we have
		\begin{equation*}
			\begin{aligned}
				\int_{\Omega'}|L(Z_{j,l})\phi|\rd x
				\le &Cd^{-(p+2)\sigma}\bigg(\sum_{i\neq j}
				\int_{\Omega_i'}\frac{1}{|x-Q_i|^{4-(p-1)\sigma+n-3\sigma}}
				\frac{1}{|x-Q_j|^{n-1}}\rd x\\
				&+\int_{\Omega_j'}\frac{1}{|x-Q_j|^{n+3-(p-1)\sigma+n-3\sigma}}
				\rd x\bigg)\\
				\le& Cd^{-n-3},
			\end{aligned}
		\end{equation*}
		and
		\begin{equation*}
			\begin{aligned}
				\int_{\omei}|L(Z_{j,l})\phi|\rd x
				\le&\frac{C}{d^{n-2}}\sum_{i\neq j}\int_{\omei_i}\frac{1}{1+|x-Q_i|^4} 	\frac{1}{1+|x-Q_j|^{n-1}}\frac{1}{1+|x-Q_i|^2}\rd x \\
				&+\frac{C}{d^{2(n-2)}}\int_{\omei_j}\frac{1}{1+|x-Q_j|^{5}}\frac{1}{1+|x-Q_j|^2}\rd x\\
				\le& \left\{\begin{aligned}
					&Cd^{-n-3}\log d, &n\ge 7,\\
					&Cd^{-8}, &n=6.
				\end{aligned}\right.
			\end{aligned}
		\end{equation*}
		Similar estimates can be obtained when $n=5$. To summarize, we have
		\begin{equation}\label{est L Zjl phi}
			\int_{\R}|L(Z_{j,l})\phi|\rd x\le
			\left\{\begin{aligned}
				&Cd^{-n-3}\log d, &n\ge 7,\\
				&Cd^{-2(n-2)}, &n=5,6.
			\end{aligned}\right.
		\end{equation}

		Similarly, using \eqref{est L Zj0}, we can obtain
		\begin{equation}\label{est L Zj0 phi}
			\int_{\R}|L(Z_{j,0})\phi|\rd x\le
			\left\{\begin{aligned}
				&Cd^{-n-2}, &n\ge 7,\\
				&Cd^{-8}\log d, &n=6,\\
				&Cd^{-6}, &n=5.
			\end{aligned}\right.
		\end{equation}

		On the other hand, following the definition of $\|\cdot\|_{**}$, it is straightforward to show that
		\begin{equation}\label{est h norm}
			\begin{split}
				\int_{\mathbb{R}^n}|h Z_{j,l}|\rd x
				\le  \frac{C}{d^{n-2}}\|h\|_{**}.
			\end{split}
		\end{equation}
		for all $n\geq 5$ and $l=0,1,2$. Hence, it follows from \eqref{est L Zjl phi}, \eqref{est L Zj0 phi} and \eqref{est h norm} that
		\begin{equation}\label{est of c}
			\|\bc\|_{\infty}\le C d^{-(n-2)}\left(d^{-1}+\|h\|_{**}\right).
		\end{equation}

		Next, we use Green's representation formula to obtain an estimate of $\phi$. Let $G(x,y)$ be the Green's function on $\mathbb{R}^n$. Then according to \eqref{line theo}, we have
		\begin{equation}\label{green function method}
			\begin{aligned}
				\phi(x)=&-p\int_{\mathbb{R}^n}\hat{K}(y)W^{p-1}(y)\phi(y) G(x,y)\rd y+
				\int_{\mathbb{R}^n}h(y)G(x,y)\rd y\\
				&+\sum_{j=1}^{k}\sum_{l=0}^{2}c_{j,l}
				\int_{\mathbb{R}^n}U_j^{p-1}(y)Z_{j,l}(y)G(x,y)\rd y.
			\end{aligned}
		\end{equation}
		We now establish the estimate on the right-hand side of \eqref{green function method} term by term.
		
		By similar argument of Lemma 3.6 in \cite{DSWduke}, we can deduce that
		\begin{equation}\label{est G1}
			\int_{\mathbb{R}^n}|h(y)G(x,y)|\rd y\le\|h\|_{**}	\int_{\mathbb{R}^n}|V(y)G(x,y)|\rd y\le C\|h\|_{**}\mathcal{W}(x).
		\end{equation}
		
		Next, taking advantage of \eqref{est of c}, we obtain
		\begin{equation}\label{est G2}
			\begin{aligned}
				\sum_{j=1}^{k}\sum_{l=0}^{2}c_{j,l}
				\int_{\mathbb{R}^n}|U_j^{p-1}Z_{j,l}G(x,y)|\rd y\le C\sum_{j=1}^k\frac{\|\bc\|_\infty}{1+|x-Q_j|^{n-2}}\le&C\left(d^{-1}+\|h\|_{**}\right)\mathcal{W}(x).
			\end{aligned}
		\end{equation}

		For the first term on the right-hand side of \eqref{green function method}, we only consider the case $n\ge 6$. Since $\hat{K}$ is bounded, using $\|\phi\|_*=1$ and estimates similar to \eqref{est W i}, we obtain 
		\begin{equation*}
			\begin{aligned}
				\big|\int_{\mathbb{R}^n}\hat{K}(y)W^{p-1}(y)\phi(y) G(x,y)\rd y\big|
				\le&C\sum_{i=1}^k\int_{\omei_i}\frac{d^{-(n-2)}}{1+|y-Q_i|^6}\frac{1}{|x-y|^{n-2}} \rd y\\
				&+C\sum_{i=1}^k\int_{\Omega_i'}\big|W^{p-1}(y)\phi(y) G(x,y)\rd y\big|.	
			\end{aligned}
		\end{equation*}
		
		The estimate of the first term on the right-hand side is easy, and we only emphasize that when estimating the second term for $x\in \Omega_j'$ we need to use \eqref{est lose} to avoid too much lose in decay. To be more precise, for $x\in \Omega_j'$, we have
		\begin{equation*}
			\begin{aligned}
				\big|\int_{\Omega'}W^{p-1}(y)\phi(y) G(x,y)\rd y\big|	\le &\sum_{i=1}^k\int_{\Omega_i'}\frac{1}{1+|y-Q_i|^{4-(p-1)\sigma_1}}\frac{d^{-2\sigma}}{1+|y-Q_i|^{n-2\sigma-\sigma_1}} \frac{1}{|x-y|^{n-2}}\rd y\\
				\le &\sum_{i=1}^k \frac{d^{-2\sigma}}{1+|x-Q_i|^{n-2}}.
			\end{aligned}
		\end{equation*}
		
		Then we can obtain
		\begin{equation}\label{est G3}
			\big|\int_{\mathbb{R}^n}\hat{K}(y)W^{p-1}(y)\phi (y)G(x,y)\rd y\big|\le C \mathcal{W}_1(x),
		\end{equation}
		where
		\begin{equation*}
			\mathcal{W}_1(x)=
			\sum_{i=1}^kd^{2-n} \frac{\log (2+|x-Q_i|)}{1+|x-Q_i|^4}\chi_{\{x\in \omei_j\}}+\sum_{j=1}^{k}\frac{d^{-2\sigma}}{1+|x-Q_j|^{n-2}}\chi_{\{x\in \Omega'\}}.
		\end{equation*}
		
		Substituting \eqref{est G1}, \eqref{est G2} and \eqref{est G3} into \eqref{green function method}, we obtain that
		\begin{equation}\label{result of phi}
			1=\|\phi\|_{*}\le C \left(\frac{\mathcal{W}_1(x)}{\mathcal{W}(x)}+
			d^{-1}+\|h\|_{**}\right).
		\end{equation}	
		Then we deduce from \eqref{result of phi} that there is $\hat{R}>0$ and some $j$ such that
		\begin{equation}\label{phi >}
			\|\phi\|_{L^\infty (B_{\hat{R}(Q_j)})}\ge a>0
		\end{equation}
		for some positive constant $a$. But $\phi(x+Q_j)$ converges (on any compact subsets of $\R$) to $\phi_{\infty}$, which is a bounded solution of the equation
		$$
		\Delta \phi_{\infty}+pU_{0,\Lambda}^{p-1}\phi_{\infty}=0.
		$$
		By passing the orthogonal conditions of $\phi$ to limit, and using the non-degeneracy of Aubin-Talenti bubble, we get that $\phi_\infty=0$. This is a contradiction to \eqref{phi >}. Thus, we deduce \eqref{priori est of phi and c} and complete the proof.
	\end{proof}

	\subsection{Nonlinear Analysis}
	We begin with the estimate of error term $E$.
	\begin{lemma}\label{lem E norm}
		Let $E$ be defined by \eqref{def E}. Then there exists a constant $C$ such that for sufficiently large $k$, we have
		\begin{equation}\label{E norm}
			\|E\|_{**}\le C.
		\end{equation}
	\end{lemma}
	\begin{proof}	
		We estimate $E$ in different regions. If $x\in \omei_i$ for some $i$, then by \eqref{est W i}, we have
		\begin{equation}\label{est E norm 1}
			\begin{split}
				|(\hat{K}-1)W^p|\le C\mu^{-m} \frac{|x-Q_i|^m}{1+|x-Q_i|^{n+2}}
				\le\frac{Cd^{-(n-2)}}{1+|x-Q_i|^4},
			\end{split}
		\end{equation}
		and
		\begin{equation}\label{est E norm 2}
			|W^p(x)-\sum_{j=1}^k U^p_{j}(x)|\le
			W^p-U_{i}^p+\frac{C}{d^{n+2}}\le CU_{i}^{p-1}\sum_{i\neq j} U_{j}+\frac{C}{d^{n+2}}
			\le \frac{Cd^{-(n-2)}}{1+|x-Q_i|^4}.
		\end{equation}
		In \eqref{est E norm 1} we also used the facts $m<n-2$ and $\mu^m\sim d^{n-2}$.
		
		If $x\in \Omega_i'$ for some $i$,  then by \eqref{est W i} again, we have
		\begin{equation*}
			|E(x)|\le 
			C(|W(x)|^p+|\sum_{j=1}^{k}U_j^{p}(x)|) \le 
			\frac{Cd^{-p\sigma}}{1+|x-Q_i|^{n+2-p\sigma}}.  
		\end{equation*}
		Thus \eqref{E norm} follows from above estimates and the definition of $\|\cdot\|_{**}$.
	\end{proof}
	
	Now we are in the position to give the proof of Proposition \ref{prop nonline proj}.
	\begin{proof}[Proof of Proposition \ref{prop nonline proj}] Recall that $\BR_*$ is defined by \eqref{def B*}. We define a map
		$$
		\mathcal{A}:=-L^{-1}(E+N(\phi)): \BR_*\to \BR_*.
		$$
		Then solving \eqref{eq phi proj} is equivalent to finding a fixed point of $\AR$. We only prove the case for $n\ge 6$, and the proof for $n=5$ only needs slight modification.
		
		For $n\ge 6$, we have $p\le 2$, so 
		$
		|N(\phi)|\le CW^{p-2}|\phi|^{2}
		$ by the mean value theorem.
		Then for $x\in \omei_i$, we have
		\begin{equation*}
			|N(\phi)|\le Cd^{4-2n}\|\phi\|_*^2 \frac{1+|x-Q_i|^{n-6}}{1+|x-Q_i|^4};
		\end{equation*}
		while for $x\in \Omega_i'$, by \eqref{est lose}, we have
		\begin{equation*}
			\begin{aligned}
				|N(\phi)|\le &
				C(1+|x-Q_i|^{n-6}) \|\phi\|_{*}^{2}\left(\sum_{i=1}^{k} \frac{d^{-2\sigma}}{1+|x-Q_i|^{n-2\sigma}}\right)^2
				\le  \frac{Cd^{-4\sigma}\|\phi\|_{*}^{2}}{1+|x-Q_i|^{n+6-4\sigma-2\sigma_1}}.
			\end{aligned}
		\end{equation*}
		Then since $\|\phi\|_{*}^{2}\le C_0$, we get
		$\|N(\phi)\|_{**}\le Cd^{-1}$.
		Therefore
		\begin{equation*}
			\|\mathcal{A}(\phi)\|_{**}\le C
			(\|E\|_{**}+\|N(\phi)\|_{**})\le C.
		\end{equation*}
		
		Next, using the fact that $|N'(t)|\le Ct^{p-1}$, we have, for $x\in \omei_j$,
		\begin{equation*}
			\begin{aligned}
				|N(\phi_1)-N(\phi_2)|
				\le&C(\|\phi_1\|_{*}^{p-1}+\|\phi_2\|_{*}^{p-1}) \|\phi_1-\phi_2\|_*(\frac{d^{2-n}}{1+|x-Q_i|^{2}})^p\\
				\le&\frac{1}{2}\|\phi_1-\phi_2\|_*(\frac{d^{2-n}}{1+|x-Q_j|^{2}}).
			\end{aligned}
		\end{equation*}
		Similarly, for $x\in \Omega_j'$, using \eqref{est lose}, we obtain
		\begin{equation*}
			\begin{aligned}
				|N(\phi_1)-N(\phi_2)|
				\le&C(\|\phi_1\|_{*}^{p-1}+\|\phi_2\|_{*}^{p-1}) \|\phi_1-\phi_2\|_*(\frac{d^{-2\sigma}}{1+|x-Q_j|^{n-2\sigma-\sigma_1}})^p\\
				\le&\frac{1}{2}\|\phi_1-\phi_2\|_*(\frac{d^{-2\sigma}}{1+|x-Q_j|^{n+2-2\sigma}})
			\end{aligned}
		\end{equation*}
		for $k$ large enough if
		$\sigma_1<2-\frac{8\sigma}{n+2}$, which can be achieved for all $n\ge 6$.
		Therefore
		\begin{equation*}
			\|\mathcal{A}(\phi_1)-\mathcal{A}(\phi_2)\|_{**}\le 
			\|N(\phi_1)-N(\phi_2)\|_{**}\le \frac{1}{2}\|\phi_1-\phi_2\|_{*}.
		\end{equation*}
		Thus, $\mathcal{A}(\phi)$ is a contraction mapping on  $\mathcal{B}_{*}$. Hence there is a unique $\phi\in \BR_*$ such that \eqref{eq phi proj} holds.
		
		By similar approach as in Proposition 4.1 in \cite{DWY}, we can prove the differentiability of $\phi(\alpha, \bq)$ as the function of $(\alpha, \bq)$  and \eqref{est of dependence}. We omit the details.
		
	\end{proof}

	\section{Computation of Projections}\label{proj}	
	~~~~This section is devoted to some preparations to the reduced problem. We define
	\begin{equation}\label{def hbc}
		\hbc=\bc-\gamma(R\bp+\bq^\perp), \text{ for every } \gamma\in\mathbb{R},
	\end{equation}
	where we write $\bp=(\mathbf{0},\mathbf{0},\mathbf{1})^{\T}\in\mathbb{R}^{3k}$ and $\bq^\perp=(0,-\bg,\bbf)$ for $\bq=(\bla,\bbf,\bg)$. Then since
	\begin{equation*}
		\frac{\partial \tw}{\partial \alpha}=(R+f_j)\frac{\partial \tw}{\partial g_j}-g_j \frac{\partial \tw}{\partial f_j}=(R\bp+\bq^\perp)\frac{\partial \tw}{\partial \bq},
	\end{equation*}
	we can write  equation \eqref{eq phi proj} in the following form
	\begin{equation}\label{eq hbc}
		L(\phi)+E+N(\phi)=\hbc\cdot \frac{\partial \tilde{W}}{\partial \bq}+\gamma \frac{\partial \tilde{W}}{\partial \alpha}.
	\end{equation}
	
	Our strategy is to solve $\hbc(\alpha,\bq,\gamma)=0$ for each $\alpha$ by adjusting $\bq$ and $\gamma$. To achieve this, we multiply \eqref{eq hbc} by $\frac{\partial {W}}{\partial \bq}$ to get
	\begin{equation*}
		\int_{\R}(	L(\phi)+E+N(\phi))\frac{\partial {W}}{\partial \bq}\rd x=M_0\hbc+\gamma\hbp,
	\end{equation*}
	where 
	\begin{equation}\label{def hbp}
		\hbp=\int_{\R}\frac{\partial \tw}{\partial \alpha} \frac{\partial W}{\partial \bq}\rd x=M_0(R\bp+\bq^{\perp}).
	\end{equation}
	Then, in view of the invertibility of $M_0$, $\hbc=0$ if and only if
	\begin{equation}\label{int hat p}
		\int_{\R}(	L(\phi)+E+N(\phi))\frac{\partial {W}}{\partial \bq}\rd x=\gamma\hbp.
	\end{equation}
	The purpose of this section is to compute each term on the left-hand side of \eqref{int hat p}.
	\subsection{Projections along $ Z_{j,1}$ and $Z_{j,2}$}
	We first compute $\int_{\R} E Z_{j,1}\rd x$ and $\int_{\R} E Z_{j,2}\rd x$ and obtain following expansions.
	\begin{lemma}\label{lem E Zj 12}
		Let $\bq\in\Xi$.  Recall that $\gamma_1$, $\gamma_2$ and $\gamma_4$ are defined by \eqref{ga 1},  \eqref{ga 2} and \eqref{ga 4} respectively. Then for sufficiently large $k$, the following expansions hold. If $n=5$, then
		\begin{equation}\label{E Zj 12}
			\begin{aligned}
				\int_{\R} E\nabla U_{j}\rd x=&
				\gamma_1\Lambda_j^{-\frac{n-2}{2}}
				\sum_{l\neq j}\Lambda_l^{-\frac{n-2}{2}}
				\frac{Q_j-Q_l}{|Q_j-Q_l|^n}+c_0\mu^{-m}\gamma_2\Lambda_j^{2-m}(|Q_j|-\mu r_0)
				\frac{Q_j}{|Q_j|}\\
				&+c_0\mu^{-m}\gamma_4\Lambda_j^{-m}\frac{Q_j}{|Q_j|^2}+O(\frac{\log d}{d^6}).
			\end{aligned}
		\end{equation}
		
		If $n\ge 6$, then for $\tau$ given by \eqref{def tau}, we have
		\begin{equation}\label{E Zj 1}
			\begin{aligned}
				\int_{\R} EZ_{j,1}\rd x=&
				\gamma_1\Lambda_j^{-\frac{n-2}{2}}
				\sum_{l\neq j}\Lambda_l^{-\frac{n-2}{2}}
				\frac{(Q_j-Q_l)\cdot\vn_j}{|Q_j-Q_l|^n}+c_0\mu^{-m}\gamma_2\Lambda_j^{2-m}(|Q_j|-\mu r_0)
				\frac{Q_j\cdot\vn_j}{|Q_j|}\\
				&+c_0\mu^{-m}\gamma_4\Lambda_j^{-m}\frac{Q_j\cdot \vn_j}{|Q_j|^2}+O(\frac{1}{d^{n+1}}+\frac{1}{d^{n+\tau}k^{2}}),
			\end{aligned}
		\end{equation}
		\begin{equation}\label{E Zj 2}
			\begin{aligned}
				\int_{\R} EZ_{j,2}\rd x=&
				\gamma_1\Lambda_j^{-\frac{n-2}{2}}
				\sum_{l\neq j}\Lambda_l^{-\frac{n-2}{2}}
				\frac{(Q_j-Q_l)\cdot\vt_j}{|Q_j-Q_l|^n}+c_0\mu^{-m}\gamma_2\Lambda_j^{2-m}(|Q_j|-\mu r_0)
				\frac{Q_j\cdot\vt_j}{|Q_j|}\\
				&+c_0\mu^{-m}\gamma_4\Lambda_j^{-m}\frac{Q_j\cdot \vt_j}{|Q_j|^2}+O(\frac{1}{d^{n+\tau}k^{2}}).
			\end{aligned}
		\end{equation}
	\end{lemma}
	
	\begin{proof}
		First we observe that the main terms in \eqref{E Zj 1} and \eqref{E Zj 2} have the same form, and can be written in one formula as in \eqref{E Zj 12}. The difference is that the small term in \eqref{E Zj 2} is smaller than that in \eqref{E Zj 1}, which is crucial for later use. Thus we only give the proof of \eqref{E Zj 2}, and the other estimates can be proved easily by applying \eqref{E1}.
		
		We will make use of symmetry repeatedly to get sharper estimate of the the small terms in \eqref{E Zj 2}. To be more precise, we define
		\begin{equation*}
			\hat{U}_j=U_{Q_{j,\alpha},\Lambda},\  \hz_{j,2}=\nabla\hu_j\cdot\vt_j,\
			\hw=\sum_{j=l}^k\hu_l,\
			\he=\hat{K}\hw^p-\sum_{l=1}^k\hu_l^p.
		\end{equation*}
		Then by symmetry,
		$\int_{\R}\he \hz_{j,2}\rd x=0$.
		Therefore $\int_{\R}EZ_{j,2}\rd x$ can be viewed as a perturbation of $0$, and we can extract extra smallness by the mean value theorem.
		
		By definition of $E$, we  write
		\begin{equation*}
			\begin{aligned}
				\int_{\R} EZ_{j,2}\rd x=&
				\int_{\R}\left(	(\hk-1)W^p\right)Z_{j,2}\rd x
				+\int_{\R} (W^p-\sum_{l=1}^{k}U_{l}^p)Z_{j,2}\rd x=:I_1+I_2.
			\end{aligned}
		\end{equation*}
		
		\noindent\textbf{Claim 1:}
		\begin{equation}\label{I1}
			I_{1}=c_0\mu^{-m}\gamma_2\Lambda_j^{2-m}(|Q_j|-\mu r_0)
			\frac{Q_j\cdot\vt_j}{|Q_j|}
			+c_0\mu^{-m}\gamma_4\Lambda_j^{-m}\frac{Q_j\cdot\vt_j}{|Q_j|^2}
			+O(\frac{1}{d^{n+\tau}k^2}).
		\end{equation}
		
		\noindent\textbf{Claim 2:}
		\begin{equation}\label{I2}
			I_2=\gamma_1\sum_{l\neq j}\Lambda_j^{-\frac{n-2}{2}}\Lambda_l^{-\frac{n-2}{2}}
			\frac{(Q_j-Q_l)\cdot\vt_j}{|Q_j-Q_l|^n}
			+O(\frac{\log d}{d^{n+2}}).
		\end{equation}
		
		The desired result follows directly from Claim 1 and Claim 2, so the remainder of this proof will be devoted to the proofs of above claims.
		
		\noindent\textbf{Proof of Claim 1:} We divide $I_1$ into three parts:
		\begin{equation*}
			\begin{aligned}
				I_1&=\int_{\omei}\left(	(\hk-1)W^p\right)Z_{j,2}\rd x
				+\int_{\omem}\left(	(\hk-1)W^p\right)Z_{j,2}\rd x
				+\int_{\omeo}\left(	(\hk-1)W^p\right)Z_{j,2}\rd x\\
				&=:I_{11}+I_{12}+I_{13}.
			\end{aligned}
		\end{equation*}

		We begin with the estimate of $I_{11}$. By $(\mathbf{K}_1)$, we can write
		\begin{equation}\label{I11 decom}
			\begin{aligned}
				I_{11}=&-c_0\mu^{-m}\sum_{i\neq j}\int_{\omei_i}||x|-\mu r_0|^m
				(W^p Z_{j,2}-\hat{W}^p\hz_{j,2})\rd x\\
				&-c_0\mu^{-m}\int_{\omei_j}||x|-\mu r_0|^m
				\left((W^p-U_j^p)Z_{j,2}-(\hw^p-\hu_j^p)\hz_{j,2}\right)\rd x\\
				&-c_0\mu^{-m}\int_{\omei_j}||x|-\mu r_0|^m
				U_j^pZ_{j,2}\rd x+O(\mu^{-m-2-\tau})\\
				=&:I_{111}+I_{112}+I_{113}+O(d^{-n-\tau}k^{-2}).
			\end{aligned}
		\end{equation}
		.
		
		By the mean value theorem and the fact that $|\bq|_\Xi\le 1$, we have
		\begin{equation}\label{W^p Z2-}
			|W^p Z_{j,2}-\hat{W}^p\hz_{j,2}|\le \left\{\begin{aligned}
				&\frac{Cd^{-p\sigma-1}}{1+|x-Q_i|^{n+2-p\sigma}} \frac{1}{1+|x-Q_j|^{n-1}}, &x\in \omem_i,\\
				&\frac{Cd^{-1}}{1+|x-Q_i|^{n+2}}\frac{1}{1+|x-Q_j|^{n-1}}, &x\in \omei_i, i\neq j.
			\end{aligned}\right.
		\end{equation}
		Then
		\begin{equation}\label{est I111}
			|I_{111}|\le C\mu^{-m}d^{-1}\sum_{i\neq j} \int_{\omei_i}\frac{||x|-\mu r_0|^m}{1+|x-Q_i|^{n+2}} \frac{1}{1+|x-Q_j|^{n-1}}\rd x\le Cd^{m-2n}.
		\end{equation}
		
		In $\omei_j$, we have
		\begin{equation*}
			\Big|(W^p-U_j^p)Z_{j,2}-(\hw^p-\hu_j^p)\hz_{j,2}\Big|\le Cd^{-1}U_j^{p-1}(\sum_{l\neq j}U_l) \frac{1}{1+|x-Q_j|^{n-1}},
		\end{equation*}
		so
		\begin{equation}\label{est I112}
			\begin{aligned}
				|I_{112}|\le C\mu^{-m}d^{-1}\sum_{i\neq j}\int_{\omei_j}\frac{|x-\mu r_0|^m}{1+|x-Q_j|^{n+3}}\frac{1}{1+|x-Q_i|^{n-2}}\rd x
				\le\left\{\begin{aligned}
					&\frac{C\log d}{d^{2n-m}}, &m\ge3,\\
					& \frac{C}{d^{2n-3}}, &m< 3.
				\end{aligned}\right. 
			\end{aligned}
		\end{equation}
		
		To estimate  $I_{113}$, we write $x=Q_{i}+y$ for $x\in \omei_i$. Then we have expansion
		\begin{equation}\label{x-mu r0}
			\begin{aligned}
				\big||x|-\mu r_0\big|^m
				=&\Big|\frac{y\cdot Q_i}{|Q_i|}\Big|^m
				+m\Big|\frac{y\cdot Q_i}{|Q_i|}\Big|^{m-2}
				\frac{y\cdot Q_i}{|Q_i|}
				\left(|Q_i|-\mu r_0
				+\frac{|y|^2}{2|Q_i|}
				-\frac{1}{8}|Q_i|(\frac{2y\cdot Q_i}{|Q_i|^2})^2
				\right)\\
				+&\Psi_i(y)+O(\frac{|y|^{m+3}}{|Q_i|^3}),
			\end{aligned}
		\end{equation}
		where $\Psi_i(y)$ is an even function in $y$.
		
		With the help of \eqref{diffenrence}, using \eqref{x-mu r0} and symmetry, we compute
		\begin{equation}\label{est I113}
			\begin{aligned}
				\int_{\omei_j}||x|-\mu r_0|^m
				U_j^pZ_{j,2}\rd x
				=&\frac{m}{p+1}(|Q_j|-\mu r_0)\frac{Q_j\cdot\vt_j}{|Q_j|}
				\int_{B_{\frac{d}{2}}}|z_1|^{m-2}z_1
				\frac{\partial U_{0,\Lambda_j}^{p+1}}{\partial z_1} \rd z+O(\frac{1}{d^{n+2-m}})\\
				&+\frac{m}{p+1}\frac{Q_j\cdot\vt_j}{|Q_j|}
				\int_{B_{\frac{d}{2}}}|z_1|^{m-2}z_1(\frac{|z|^2}{2|Q_j|}-\frac{z_1^2}{2|Q_j|})
				\frac{\partial U_{0,\Lambda_j}^{p+1}}{\partial z_1}\rd z+O(\frac{1}{\mu^3})\\
				=&-\frac{m(m-1)}{p+1}\Lambda_j^{2-m}(|Q_j|-\mu r_0)\frac{Q_j\cdot\vt_j}{|Q_j|}
				\int_{\R}|z_1|^{m-2} U^{p+1}(z)\rd z+O(\frac{1}{d^{n+2-m}})\\
				&-\frac{m(m-1)}{p+1}\Lambda_j^{-m}\frac{Q_j\cdot\vt_j}{2|Q_j|^2}
				\int_{\R}|z_1|^{m-2}(|z|^2-z_1^2)
				U^{p+1}(z)\rd z+O(\frac{1}{\mu^3}).
			\end{aligned}
		\end{equation}
		
		In conclusion, combining \eqref{est I111}, \eqref{est I112} and \eqref{est I113}, we obtain
		\begin{equation}\label{est I11}
			I_{11}=c_0\mu^{-m}\gamma_2\Lambda_j^{2-m}(|Q_j|-\mu r_0)
			\frac{Q_j\cdot\vt_j}{|Q_j|}
			+c_0\mu^{-m}\gamma_4\Lambda_j^{-m}\frac{Q_j\cdot\vt_j}{|Q_j|^2}
			+O(\frac{1}{d^{n+\tau}k^2})
		\end{equation}
		
		Next we consider $I_{12}$. As in \eqref{I11 decom},  we can write
		\begin{equation*}
			\begin{aligned}
				I_{12}=&-c_0\mu^{-m}\sum_{i=1}^k\int_{\omem_i}||x|-\mu r_0|^m
				(W^p Z_{j,2}-\hat{W}^p\hz_{j,2})\rd x+O(d^{-n-\tau}k^{-2}).
			\end{aligned}
		\end{equation*}
		Then using \eqref{W^p Z2-} again, we get that
		\begin{equation}\label{est I12}
			\begin{aligned}
				|I_{12}|\le& C\mu^{-m}d^{-p\sigma-1}\sum_{i=1}^k \int_{\omem_i}\frac{||x|-\mu r_0|^m}{1+|x-Q_i|^{n+2-p\sigma}} \frac{1}{1+|x-Q_j|^{n-1}}\rd x+Cd^{-n-\tau}k^{-2}\\
				\le& Cd^{m-2n}+Cd^{-n-\tau}k^{-2}.
			\end{aligned}
		\end{equation}
		
		Finally for $I_{13}$, we have
		\begin{equation}\label{est I13}
			\begin{aligned}
				|I_{13}|&\le \sum_{i\neq j} \frac{C}{|Q_i-Q_j|^{\sigma}} 
				\int_{\omeo_i}\frac{d^{-p\sigma}}{1+|x-Q_i|^{2n+1-(p+1)\sigma}}
				\rd x+\frac{Cd^{-p\sigma}}{\mu^{n+1-p\sigma}}\le \frac{C}{d^{n+1}k^{n+1-(p+1)\sigma}}.
			\end{aligned}
		\end{equation}	
		Then \eqref{I1} follows from \eqref{est I11}, \eqref{est I12} and \eqref{est I13}.\\

		\noindent\textbf{Proof of Claim 2:}
		We also divide $I_2$ into three parts:
		\begin{equation*}
			\begin{aligned}
				I_2=&\int_{\omei_j} pU_{j}^{p-1}\sum_{l\neq j}U_{l}Z_{j,2}\rd x
				+\int_{\omei_j} (W^p-\sum_{l=1}^{k}U_{l}^p-pU_{j}^{p-1}\sum_{l\neq j}U_{l})Z_{j,2}\rd x
				+\int_{\Omega_j'} (W^p-\sum_{l=1}^{k}U_{l}^p)Z_{j,2}\rd x\\
				=&:I_{21}+I_{22}+I_{23}.
			\end{aligned}
		\end{equation*}
		
		First we consider $I_{21}$. For simplicity, we denote $\vec\nu_{jl}=\frac{Q_j-Q_l}{|Q_j-Q_l|}$. Then by a change of variable and utilizing symmetry, we have
		\begin{equation*}
			\begin{aligned}
				\int_{\omei_j} \frac{\partial U_j^p}{\partial \vt_j}U_{l}\rd x=&(\vec\nu_{jl}\cdot\vt_j)
				\int_{B_{\frac{d}{2}}}
				\frac{\partial U_{0,\Lambda_j}^p(y)}{\partial y_1}U_{0,\Lambda_l}(y-\rho_{jl}\vec e_1)\rd y+O(\frac{1}{d_{jl}^{n-2}}\frac{1}{d^4})\\
				=&(\vec\nu_{jl}\cdot\vt_j)\left(\int_{B_{\frac{d}{2}}}
				\frac{\partial U_{0,\Lambda_j}^p(y)}{\partial y_1}y_1\frac{\partial U_{0,\Lambda_j}(\rho_{jl}\vec e_1)}{\partial y_1}\rd y-\int_{\R\setminus B_{\frac{d}{2}}}
				\frac{\partial U_{0,\Lambda_j}^p(y)}{\partial y_1}y_1\frac{\partial U_{0,\Lambda_j}(\rho_{jl}\vec e_1)}{\partial y_1}\rd y\right)\\
				&+(\vec\nu_{jl}\cdot\vt_j)\int_{B_{\frac{d}{2}}}
				\frac{\partial U_{0,\Lambda_j}^p(y)}{\partial y_1}\Bigg(U_{0,\Lambda_l}(y-\rho_{jl}\vec e_1)-U_{0,\Lambda_l}(\rho_{jl})+y\cdot \nabla U_{0,\Lambda_l}(\rho_{jl}\vec e_1)\\
				&\quad\quad-\frac{1}{2}y^{\T}D^2U_{0,\Lambda_l}(\rho_{jl}\vec e_1)y\Bigg)\rd y
				+O(\frac{1}{d_{jl}^{n-2}}\frac{1}{d^4}).
			\end{aligned}			
		\end{equation*}
		By the mean value theorem,	we have
		\begin{equation*}			\begin{aligned}
				&(\vec\nu_{jl}\cdot\vt_j)\frac{\partial U_{0,\Lambda_j}(\rho_{jl}\vec e_1)}{\partial y_1}\int_{\R\setminus B_{\frac{d}{2}}}
				y_1\frac{\partial U_{0,\Lambda_j}^p(y)}{\partial y_1}\rd y\\
				=&	\cos\frac{\theta_j-\theta_l}{2}
				\frac{\partial U_{0,\Lambda}(d_{jl}\vec e_1)}{\partial y_1}\int_{\R\setminus B_{\frac{d}{2}}}
				y_1\frac{\partial U_{0,\Lambda}^p(y)}{\partial y_1}\rd y+O(d^{-n-2}),
		\end{aligned}		\end{equation*}
		and
		\begin{equation*}
			\begin{aligned}
				&(\vec\nu_{jl}\cdot\vt_j)\int_{B_{\frac{d}{2}}}
				\frac{\partial U_{0,\Lambda_j}^p(y)}{\partial y_1}\Bigg(U_{0,\Lambda_l}(y-\rho_{jl}\vec e_1)-U_{0,\Lambda_l}(\rho_{jl})+y\cdot \nabla U_{0,\Lambda_l}(\rho_{jl}\vec e_1)-\frac{1}{2}y^{\T}D^2U_{0,\Lambda_l}(\rho_{jl}\vec e_1)y\Bigg)\rd y\\
				=&\cos\frac{\theta_j-\theta_l}{2}
				\int_{B_{\frac{d}{2}}}\frac{\partial U_{0,\Lambda}^p(y)}{\partial y_1}\Bigg(U_{0,\Lambda_l}(y-\rho_{jl}\vec e_1)-U_{0,\Lambda_l}(\rho_{jl})+y\cdot \nabla U_{0,\Lambda_l}(\rho_{jl}\vec e_1)-\frac{1}{2}y^{\T}D^2U_{0,\Lambda_l}(\rho_{jl}\vec e_1)y\Bigg)\rd y\\
				&+ O(d_{jl}^{-n-2}\log d +d_{jl}\mu^{-1}d^{-n-2}\log d).
			\end{aligned}
		\end{equation*}
		
		Taking summation of above estimates in $l$ and using symmetry, we obtain	
		\begin{equation}\label{est I21}
			\begin{aligned}
				I_{21}=\sum_{l\neq j}\int_{\omei_j} \frac{\partial U_j^p}{\partial \vt_j}U_{l}\rd x
				=&\sum_{l\neq j}(\vec\nu_{jl}\cdot\vt_j)\frac{\partial U_{0,\Lambda_j}(\rho_{jl}\vec e_1)}{\partial y_1}\int_{\R}
				y_1\frac{\partial U_{0,\Lambda_j}^p(y)}{\partial y_1}\rd y+O(\frac{\log d}{d^{n+2}})\\
				=&\gamma_1\sum_{l\neq j}\Lambda_j^{-\frac{n-2}{2}}\Lambda_l^{-\frac{n-2}{2}}
				\frac{(Q_j-Q_l)\cdot \vt_j}{|Q_j-Q_l|^n}
				+O(\frac{\log d}{d^{n+2}})
			\end{aligned}
		\end{equation}

		For $I_{22}$, we observe that in $\omei_j$ there holds
		\begin{equation*}
			\begin{aligned}
				\Big|(W^p-\sum_{l=1}^{k}U_{l}^p-pU_{j}^{p-1}\sum_{l\neq j}U_{l})Z_{j,2}
				-(\hw^p-\sum_{l=1}^{k}\hu_{l}^p-p\hu_{j}^{p-1}\sum_{l\neq j}\hu_{l})\hz_{j,2}\Big|
				\le\frac{Cd^{2(2-n)-1}}{1+|x-Q_j|^{5}}.
			\end{aligned}
		\end{equation*}
		Then by symmetry, we have
		\begin{equation}\label{est I22}
			\begin{aligned}
				I_{22}=O(\int_{\omei_j}\frac{d^{2(2-n)-1}}{1+|x-Q_j|^{5}})
				= O(\frac{1}{d^{n+2}}).
			\end{aligned}
		\end{equation}
		
		For $I_{23}$, we notice that similar to \eqref{W^p Z2-}, there holds
		\begin{equation*}
			|(W^{p}-\sum_{l=1}^kU_l^p) Z_{j,2}-(\hat{W}^p-\sum_{l=1}^k\hu_l^p)\hz_{j,2}|\le \left\{\begin{aligned}
				&\frac{Cd^{-p\sigma-1}}{1+|x-Q_i|^{n+2-p\sigma}} \frac{1}{1+|x-Q_j|^{n-1}}, &x\in \Omega'_i,\\
				&\frac{Cd^{1-n}}{1+|x-Q_i|^{4}}\frac{1}{1+|x-Q_j|^{n-1}}, &x\in \omei_i, i\neq j.
			\end{aligned}\right.
		\end{equation*}
		Then by symmetry again,
		\begin{equation}\label{est I23}
			\begin{aligned}
				|I_{23}|
				\le&\sum_{i=1}^k\int_{\Omega_i'}\frac{Cd^{-p\sigma-1}}{1+|x-Q_i|^{n+2-p\sigma}} \frac{1}{1+|x-Q_j|^{n-1}}\rd x+\sum_{i\neq j}\int_{\omei_i}\frac{Cd^{1-n}}{1+|x-Q_i|^{4}}\frac{1}{1+|x-Q_j|^{n-1}}\rd x\\
				\le& \frac{C}{d^{n+2}}.
			\end{aligned}
		\end{equation}
		Then \eqref{I2} follows from \eqref{est I21}, \eqref{est I22} and \eqref{est I23}.

	\end{proof}
	
	Next we compute $\int_{\R} (L(\phi)+N(\phi))Z_{j,l}\rd x$ for $l=1,2$.
	\begin{lemma}\label{lem L+N Zj 12}
		Let $l=1,2$. Then for sufficiently large $k$, we have
		\begin{equation}\label{proj L+N Zjl}
			|\int_{\R} (L(\phi)+N)Z_{j,l}\rd x|\le
			\left\{\begin{aligned}
				&Cd^{-n-2}, &n\ge 6,\\
				&Cd^{-6}, &n=5.
			\end{aligned}\right.
		\end{equation}
	\end{lemma}
	\begin{proof}
		Since $\|\phi\|_*\le C$, the estimate of $\int_{\R} L(\phi)Z_{j,l}\rd x$ is already proved in \eqref{est L Zjl phi}.
		
		For $\int_{\R} N(\phi)Z_{j,l}\rd x$, we observe that if $n\ge 6$, then $|N(\phi)|\le CW^{p-2}(x)|\phi|^{2}$, so
		\begin{equation*}
			\begin{aligned}
				\int_{\omei} |N(\phi)Z_{j,l}|\rd x\le \frac{C}{d^{2n-4}}\sum_{i=1}^k \int_{\omei_i}
				\frac{|x-Q_i|^{n-6}}{1+|x-Q_i|^{4}}\frac{1}{1+|x-Q_j|^{n-1}}\rd x\le \frac{C}{d^{n+5}},
			\end{aligned} 
		\end{equation*}
		and by \eqref{est W i},
		\begin{equation*}
			\begin{aligned}
				\int_{\Omega'} |N(\phi)Z_{j,l}|\rd x
				&\le C\sum_{i=1}^k \int_{\Omega_i'} \frac{1}{1+|x-Q_i|^{6-n}}	\left(\sum_{l=1}^k\frac{d^{-2\sigma}}{|x-Q_l|^{n-2\sigma}}\right)^2\frac{1}{|x-Q_j|^{n-1}}\rd x\\
				&\le \frac{C}{d^{n+5}}.
			\end{aligned} 
		\end{equation*}
		The case $n=5$ can be dealt with similarly and we omit the details.
	\end{proof}

	\subsection{Projection along $Z_{j,0}$}
	In this subsection, we compute the projection along $Z_{j,0}$. We start with $\int_{\R}EZ_{j,0}\rd x$.
	\begin{lemma}\label{lem E Zj 0}
		Let $\bq\in\Xi$. Recall that $\gamma_1$ and $\gamma_3$ are defined by \eqref{ga 1}  and \eqref{ga 3} respectively. Then for sufficiently large $k$, we have the following expansions. If $n=5$, then
		\begin{equation}\label{E Zj 0 5}
			\begin{aligned}
				\int_{\R} EZ_{j,0}\rd x=&
				-\frac{1}{2}\gamma_1\Lambda_j^{-\frac{n}{2}}
				\sum_{l\neq j}\Lambda_l^{-\frac{n-2}{2}}
				\frac{1}{|Q_j-Q_l|^{n-2}}+c_0\mu^{-m}\gamma_3\Lambda_j^{-1-m}+O(\frac{1}{d^{n}}).
			\end{aligned}
		\end{equation}
		If $n\ge 6$, then	\begin{equation}\label{E Zj 0}
			\begin{aligned}
				\int_{\R} EZ_{j,0}\rd x=&
				-\frac{1}{2}\gamma_1\Lambda_j^{-\frac{n}{2}}
				\sum_{l\neq j}\Lambda_l^{-\frac{n-2}{2}}
				\frac{1}{|Q_j-Q_l|^{n-2}}+c_0\mu^{-m}\gamma_3\Lambda_j^{-1-m}+\hat{I}+O(\frac{\log d}{d^{n+1}}).
			\end{aligned}
		\end{equation}
		Here $\hat I$ is a quantity that is independent of $j$ and satisfies $|\hat I|\le \frac{C\log d}{d^n}$.
	\end{lemma}
	\begin{proof}
		As before, we decompose
		\begin{equation*}
			\begin{aligned}
				\int_{\R} EZ_{0,j}\rd x=
				\int_{\R}	(\hk-1)W^p Z_{0,j}\rd x
				+\int_{\R} (W^p-\sum_{l=1}^{k}U_{l}^p)Z_{0,j}\rd x=:\tilde I_1+\tilde{I}_2.
			\end{aligned}
		\end{equation*}
		For the same reason as in Lemma \ref{lem E Zj 12}, we only give the proof of \eqref{E Zj 0}, which will be achieved by proving the following claims:
		
		\noindent\textbf{Claim 3:}
		\begin{equation}\label{ti1}
			\ti_1=c_0\mu^{-m}\gamma_3\Lambda_j^{-1-m}+\hat{I}_1+O(\frac{1}{d^{n+1}}).
		\end{equation}
		\noindent\textbf{Claim 4:}
		\begin{equation}\label{ti2}
			\tilde{I}_{2}=\frac{1}{2}\gamma_1
			\Lambda_j^{-\frac{n}{2}}\sum_{l\neq j}\Lambda_l^{-\frac{n-2}{2}}|Q_j-Q_l|^{2-n}+\hat I_{2}+O( d^{-n-1}\log d).
		\end{equation}
		
		\noindent\textbf{Proof of Claim 3:} We divide $\ti_1$ into three parts:
		\begin{equation*}
			\begin{aligned}
				\ti_1&=\int_{\omei}(\hk-1)W^p\rd x
				+\int_{\omem}(\hk-1)W^pZ_{j,0}\rd x
				+\int_{\omeo}(\hk-1)W^pZ_{j,0}\rd x=:\ti_{11}+\ti_{12}+\ti_{13}.
			\end{aligned}
		\end{equation*}
		We first consider $\ti_{11}$. Similar to \eqref{I11 decom}, we can write
		\begin{equation}\label{ti11 decom}
			\begin{aligned}
				I_{11}=&-c_0\mu^{-m}\sum_{i\neq j}\int_{\omei_i}||x|-\mu r_0|^m
				(W^p Z_{j,0}-\hat{W}^p\hz_{j,0})\rd x\\
				&-c_0\mu^{-m}\int_{\omei_j}||x|-\mu r_0|^m
				\left((W^p-U_j^p)Z_{j,0}-(\hw^p-\hu_j^p)\hz_{j,0}\right)\rd x\\
				&-c_0\mu^{-m}\int_{\omei_j}||x|-\mu r_0|^m
				U_j^pZ_{j,0}\rd x+O(d^{-n-\tau}k^{-2}).
			\end{aligned}
		\end{equation}
		
		Similar to \eqref{est I111} and \eqref{est I112}, we have
		\begin{equation*}
			\begin{aligned}
				\sum_{i\neq j}\int_{\omei_i}||x|-\mu r_0|^m
				(W^p Z_{j,0}-\hat{W}^p\hz_{j,0})\rd x\le& \sum_{i\neq j}\int_{\omei_i}\frac{Cd^{-1}}{1+|x-Q_i|^{n-m+2}} \frac{1}{|x-Q_j|^{n-2}}\rd x\le\frac{C}{d^{n+1-m}},
			\end{aligned}
		\end{equation*}
		and
		\begin{equation*}
			\int_{\omei_j}||x|-\mu r_0|^m
			\left((W^p-U_j^p)Z_{j,0}-(\hw^p-\hu_j^p)\hz_{j,0}\right)\rd x=O(\frac{1}{d^{n+1-m}}).
		\end{equation*}
		Using expansion in \eqref{x-mu r0}, we can obtain
		\begin{equation*}
			\begin{aligned}
				\int_{\omei_j}\left(||x|-\mu r_0|^m
				U_{j}^p\right) Z_{j,0}\rd x=&
				\int_{B_{\frac{d}{2}}}\left(||y+Q_j|-\mu r_0|^m
				U_{0,\Lambda_j}^p\right) Z_{j,0}(y+Q_j)\rd y+O(\frac{1}{d^{n+1-m}})\\
				=&\frac{1}{p+1}\int_{B_{\frac{d}{2}}} |y_1|^m\frac{\partial U_{0,\Lambda_j}^{p+1}(y)}{\partial \Lambda_j}\rd y+O((|Q_i|-\mu r_0)^2)+O(\frac{1}{d^{n+1-m}})\\
				=&-\Lambda_j^{-1-m}\gamma_3-\frac{1}{p+1}\int_{\R\setminus B_{\frac{d}{2}}}|y_1|^m\frac{\partial U_{0,\Lambda}^{p+1}(y)}{\partial \Lambda}\rd y+O(\frac{1}{\mu^2})+O(\frac{1}{d^{n+1-m}}).
			\end{aligned}
		\end{equation*}
		
		Combining above estimates, we obtain
		\begin{equation}\label{est ti11}
			\tilde{I}_{11}=c_0\mu^{-m}\Lambda_j^{-1-m}\gamma_3+\hat{I}_{11}+O(\frac{1}{d^{n+\tau}k^2}),
		\end{equation}
		where
		\begin{equation*}
			\begin{aligned}
				\hat{I}_{11}=&-c_0\mu^{-m}\sum_{i\neq j}\int_{\omei_i}||x|-\mu r_0|^m\hat{W}^p\hz_{j,0}\rd x
				-c_0\mu^{-m}\int_{\omei_j}||x|-\mu r_0|^m
				(\hw^p-\hu_j^p)\hz_{j,0}\rd x\\
				&+\frac{1}{p+1}c_0\mu^{-m}\int_{\R\setminus B_{\frac{d}{2}}}|y_1|^m\frac{\partial U_{0,\Lambda}^{p+1}(y)}{\partial \Lambda}\rd y.
			\end{aligned}
		\end{equation*}
		Similar to \eqref{est I12} and \eqref{est I13}, we can show that
		\begin{equation}\label{est ti12}
			|\tilde{I}_{12}|=\int_{\omem}||x|-\mu r_0|^m	\hw^p \hz_{j,0}\rd x+O(\frac{1}{d^{2n-m-1}})+O(\frac{1}{d^{n+\tau}k^2}),
		\end{equation}
		and
		\begin{equation}\label{est ti13}
			\begin{aligned}
				|\ti_{13}|=O(\frac{1}{d^{n}k^{n-(p+1)\sigma}}).
			\end{aligned}
		\end{equation} 
		Let $\hat{I}_1=\hat{I}_{11}+\int_{\omem}||x|-\mu r_0|^m	\hw^p \hz_{j,0}\rd x$. Then it is easy to see that $\hat{I}_1$ is independent of $j$ and  $|\hat{I}_1|\le\frac{C}{d^n}$. Thus \eqref{ti1} follows from \eqref{est ti11}, \eqref{est ti12} and \eqref{est ti13}.\\
		
		\noindent\textbf{Proof of Claim 4:} We also divide $\tilde{I}_2$ into three parts:
		\begin{equation*}
			\begin{aligned}
				\tilde I_2=&\int_{\omei_j} pU_{j}^{p-1}\sum_{l\neq j}U_{l}Z_{j,0}\rd x
				+\int_{\omei_j} (W^p-\sum_{l=1}^{k}U_{l}^p-pU_{j}^{p-1}\sum_{l\neq j}U_{l})Z_{j,0}\rd x
				+\int_{\Omega_j'} (W^p-\sum_{l=1}^{k}U_{l}^p)Z_{j,0}\rd x\\
				=&:\tilde I_{21}+\tilde I_{22}+\tilde I_{23}.
			\end{aligned}
		\end{equation*}
		
		To estimate $\tilde I_{21}$, by symmetry and \eqref{diffenrence}, we have
		\begin{equation*}
			\begin{aligned}
				\ti_{21}=&
				\sum_{l\neq j}\int_{B_{\frac{d}{2}}}
				\frac{\partial U_{0,\Lambda_j}^p(y)}{\partial \Lambda_j}\left(U_{0,\Lambda_l}(y-(Q_l-Q_j)))-U_{0,\Lambda_l}(\rho_{jl})+\nabla U_{0,\Lambda_l}(Q_l-Q_j)\cdot y\right)\rd y\\
				&+\sum_{l\neq j}\left(U_{0,\Lambda_l}(\rho_{jl})\int_{\R}
				\frac{\partial U_{0,\Lambda_j}^p(y)}{\partial \Lambda_j}\rd y-U_{0,\Lambda_l}(\rho_{jl})\int_{\R\setminus B_{\frac{d}{2}}}
				\frac{\partial U_{0,\Lambda_j}^p(y)}{\partial \Lambda_j}\rd y\right)+O(\frac{1}{d^{n+1}}).
			\end{aligned}
		\end{equation*}
		Let
		\begin{equation*}
			\begin{aligned}
				\hat{I}_{21}=&\sum_{l\neq j}\int_{B_{\frac{d}{2}}}
				\frac{\partial U_{0,\Lambda}^p(y)}{\partial \Lambda}\left(U_{0,\Lambda}(y-\hat{P}_{jl}))-U_{0,\Lambda}(d_{jl})+\nabla U_{0,\Lambda}(\hat{P}_{jl})\cdot y\right)\rd y\\
				&+\int_{\R\setminus B_{\frac{d}{2}}}
				\frac{\partial U_{0,\Lambda}^p(y)}{\partial \Lambda}\rd y \sum_{l\neq j}U_{0,\Lambda}(d_{jl}).
			\end{aligned}
		\end{equation*}
		Then $\hat{I}_{21}$ is independent of $j$ and
		$|\hat{I}_{21}|\le C\log d\sum_{l\neq j}d_{jl}^{-n}+Cd^{-2}\sum_{l\neq j}d_{jl}^{2-n}\le Cd^{-n}\log d.$
		Thus, similar to \eqref{est I21}, by the mean value theorem, we can obtain
		\begin{equation}\label{est ti21}
			\begin{aligned}
				\tilde{I}_{21}=&\sum_{l\neq j}U_{0,\Lambda_l}(\rho_{jl})\int_{\R}
				\frac{\partial U_{0,\Lambda_j}^p(y)}{\partial \Lambda_j}\rd y+\hat I_{21}+O( d^{-n-1}\log d)\\
				=&\frac{1}{2}\gamma_1
				\Lambda_j^{-\frac{n}{2}}\sum_{l\neq j}\Lambda_l^{-\frac{n-2}{2}}|Q_j-Q_l|^{2-n}+\hat I_{21}+O( d^{-n-1}\log d).    
			\end{aligned}
		\end{equation}

		Similar to estimates \eqref{est I22} and \eqref{est I23}, we have
		\begin{equation}\label{est ti22 23}
			\tilde I_{22}=\hat{I}_{22}+O(\frac{1}{d^{n+1}}), \quad \tilde I_{23}=\hat{I}_{23}+O(\frac{1}{d^{n+1}}),
		\end{equation}
		where
		\begin{equation*}
			\hat{I}_{22}=\int_{\omei_j}(\hw^p-\sum_{l=1}^{k}\hu_{l}^p-p\hu_{j}^{p-1}\sum_{l\neq j}\hu_{l})\hz_{j,0}\rd x, \ \hat{I}_{23}=\int_{\Omega_j'}(\hw^p-\sum_{l=1}^k\hu_l^p)\hz_{j,0}\rd x
		\end{equation*}
		are both independent of $j$, and it is standard to show that $|\hat{I}_{22}|+|\hat{I}_{23}|\le Cd^{-n}$.
		
		By letting $\hat{I}_{2}=\hat{I}_{21}+\hat{I}_{22}+\hat{I}_{23}$, \eqref{ti2} then follows from \eqref{est ti21} and \eqref{est ti22 23}.
	\end{proof}

	Next, for $\int_{\R} (L(\phi)+N(\phi))Z_{j,0}\rd x$, we have the following lemma, whose proof is similar to that of Lemma \ref{lem L+N Zj 12} and is omitted.
	\begin{lemma}\label{lem L+N Zj 0}
		For sufficiently large $k$, we have
		\begin{equation}\label{proj L+N Zj0}
			\Big|\int_{\R} (L(\phi)+N(\phi))Z_{j,0}\rd x\Big|\le
			\left\{\begin{aligned}
				&Cd^{-n-2}\log d, &n\ge 6,\\
				&Cd^{-6}, &n=5.
			\end{aligned}\right.
		\end{equation}
	\end{lemma}

	\section{Proof of Theorem \ref{thm main}}\label{pf of thm main}
	\subsection{Reduction to one dimension}\label{redu 1 dim}
	~~~~In this subsection, we are going to solve $\hbc(\alpha,\bq,\gamma)=0$ in \eqref{eq hbc}. This is done by a careful analysis of the invertibility of the matrix $T$ defined by \eqref{def T} below.
	
	Using results in the previous sections, we first derive an expansion of the left-hand side of \eqref{int hat p}. The proof is tedious so we leave it to Appendix \ref{pf lem proj exp}.
	\begin{lemma}\label{lem proj exp}
		For sufficiently large $k$, if $n\ge 6$, then we have the following expansion:
		\begin{equation}\label{proj exp}
			\begin{aligned}
				\int_{\R} (E+L(\phi)+N(\phi))\frac{\partial W}{\partial\bq}\rd x=
				-\frac{n-2}{2}\gamma_1\Lambda^{2-n}(2R)^{1-n}T\textbf{q}
				+\begin{pmatrix}
					\hat{I}+O(d^{-n-\frac{3}{2}\tau_1}k^{-2}+d^{-n-1}{\log d})\\
					O(d^{1-n-2\tau_1}k^{-1})\\
					O(d^{-n-2\tau_1}k^{-2})
				\end{pmatrix},
			\end{aligned}
		\end{equation}
		where $\hat{I}$ is described in Lemma \ref{lem E Zj 0}, and the $3k\times 3k$ matrix $T$ is defined by
		\begin{equation}\label{def T}
			T=\begin{pmatrix} 
				c_1A_1+c_2I & O &  A_2 \\
				O & c_4I & O \\
				\frac{1}{2} A_2 & O & c_3A_3 \\
			\end{pmatrix}.
		\end{equation}
		Here $I$ is the $k\times k$ identity matrix. $A_1, A_2,  A_3$  are $k\times k$ circulant matrices, and their entries in the first row are given by
		\begin{equation*}
			A_{1,1l}=\frac{1}{\sin^{n-2}((l-1)\theta)},\ 
			A_{2,1l}= \frac{\sin(2(l-1)\theta)}{\sin^{n}((l-1)\theta)},\ 
			A_{3,1l}=\frac{1-\frac{n-2}{n-1}\sin^{2}((l-1)\theta)}{\sin^{n}((l-1)\theta)}, l=2,\cdots,k,
		\end{equation*}
		where we recall $\theta=\frac{2\pi}{k}$, and
		\begin{equation*}
			A_{1,11}=A_{2,11}=0,\  A_{3,11}=-\sum_{l=2}^k A_{3,1l}.
		\end{equation*}
		$c_1, c_2, c_3, c_4$ are constants given by
		\begin{equation*}
			c_1={2R},\ c_2=\frac{n-2-2m}{n-2}(2R)\sum_{l=2}^k \frac{1}{\sin^{n-2}((l-1)\theta)},\
			c_3=\frac{n-1}{n-2}R^{-1},\ c_4=\frac{2c_0\mu^{-m}\gamma_2\Lambda^{{n}-m}}{(n-2)\gamma_1(2R)^{1-n}}.
		\end{equation*}
		When $n=5$, the remainder term in \eqref{proj exp} is replaced by
		\begin{equation*}
			(O(d^{-6}\log d),\ O(d^{-4-2\tau_1}k^{-1}), \ O(d^{-5}))^\T.
		\end{equation*}
	\end{lemma}
	
	Our next goal is to study the invertibility of the matrix $T$. We observe that the matrix $T$ can be decomposed into two blocks, so we first study the linear problem \eqref{eq T'v} below and establish the following lemma, whose proof is delayed to Appendix \ref{pf lem sol T'}.
	\begin{lemma}\label{lem sol T'}
		Let $\bp'=({\bf 0},{\bf1})^{\textnormal{T}}\in {\mathbb{R}}^{2k}$ and \begin{equation*}
			T'=\begin{pmatrix} 
				c_1A_1+c_2I&  A_2 \\
				\frac{1}{2} A_2 &  c_3A_3 \\
			\end{pmatrix}.
		\end{equation*}
		Then for any ${\bf{b}}=({\bf{b}}_0, {\bf{b}}_2)\in  {\mathbb{R}}^{2k}$, there exist a unique vector ${\bf v}=({\bf{v}}_0,{\bf{v}}_2)\in {\mathbb{R}}^{2k}$ and a
		unique constant $\gamma\in\mathbb{R}$	such that
		\begin{equation}\label{eq T'v}
			T'{\bf v}={\bf b}+\gamma \bp',\quad {\bf v}\perp \bp'.
		\end{equation}	
		Moreover, we have estimates
		\begin{equation}\label{est v0}
			\|{\bf v}_0\|_{\infty}\le Cd^{-1}k^{1-n}\|{\bf b}_0\|_{\infty}+Cd^{-1}k^{2-n}\|{\bf b}_0\|_{\infty'}+Ck^{2-n}\log k\|{\bf b}_2\|_{\infty},
		\end{equation}
		\begin{equation}\label{est v2}
			\|{\bf v}_2\|_{\infty}\le Ck^{2-n}\log k\|{\bf b}_0\|_{\infty}+Cdk^{3-n}\|{\bf b}_2\|_{\infty},
		\end{equation}
		and for $j\neq t$,
		\begin{equation}\label{est v2-}
			\begin{aligned}
				|v_{2,j}-v_{2,t}|\le& Ck^{2-n}\log(2+\min\{|j-t|,k-|j-t|\})\|{\bf b}_0\|_{\infty'}\\
				&+Cdk^{2-n}\min\{|j-t|,k-|j-t|\}\|{\bf b}_2\|_{\infty}.
			\end{aligned}
		\end{equation}
	\end{lemma}
	
	Using Lemma \ref{lem sol T'}, we can obtain the following result concerning the invertibility of $T$.
	\begin{proposition}\label{prop sol T}
		Let $\bp=(\mathbf{0},\mathbf{0},\mathbf{1})^{\T}\in\mathbb{R}^{3k}$ and $\hbp$ be defined by \eqref{def hbp}. Then for sufficiently large $k$ and any ${\bf b}=({\bf b}_0, {\bf b}_1, {\bf b}_2)\in \mathbb{R}^{3k}$, there exist a unique vector ${\bf v}=({\bf v}_0, {\bf v}_1, {\bf v}_2)\in \mathbb{R}^{3k}$ and a unique constant $\gamma\in \mathbb{R}$ such that
		\begin{equation}\label{eq Tv q1}
			T{\bf v}={\bf b}+\gamma \hbp, \quad {\bf v}\perp \bp.
		\end{equation}
		Moreover, estimates \eqref{est v0}-\eqref{est v2-} holds, and we also have 
		\begin{equation}\label{est v1 new}
			\|{\bf v}_1\|_{\infty}\le Cd^{-1}k^{1-n}(\|{\bf b}_1\|_{\infty}+d^{-1}k^{-1}\|{\bf b}_2\|_{\infty}),
		\end{equation}
	\end{proposition}
	\begin{proof}
		Firstly, from Lemma \ref{lem sol T'}, we know that the equation
		\begin{equation*}
			T{\bf v}={\bf b}+\gamma_0 {\bp}, \quad {\bf v}\perp {\bp}
		\end{equation*}
		has a unique solution that satisfies \eqref{est v0}-\eqref{est v2-} and
		$\|{\bf v}_1\|_{\infty}\le Cd^{-1}k^{1-n}\|{\bf b}_1\|_{\infty}$.
		
		Now we go back to equation \eqref{eq Tv q1}. The existence is obvious, and it suffices to prove the estimates. By orthogonality it is easy to see that $\gamma=-\frac{{\bf b}\cdot{\bp}}{\hbp\cdot \bp}$. By definition of $\hbp$ and \eqref{almost ortho}, we have $\hbp\cdot \bp\ge CRk$, hence  $|\gamma|\le CR^{-1}\|{\bf b}_2\|_{\infty}$.
		Thus the estimates follows easily.
	\end{proof}
	
	Since $\hbp$ depends on $\bq$, we shall denote the inverse of $T$ in the sense of Proposition \ref{prop sol T} by $T_{\bq}^{-1}$. The next corollary characterize this dependence in more detail.
	\begin{corollary}\label{coro Tq -1}
		Let $\bq, \mathring{\bq}\in \Xi$. Then for fixed ${\bf b}=({\bf b}_0, {\bf b}_1, {\bf b}_2)\in \mathbb{R}^{3k}$, we have
		\begin{equation*}
			\|(T^{-1}_{\bq}-T^{-1}_{	\mathring{\bq}}){\bf b}\|_\Xi\le  Cd^{-1+\frac{1}{2}\tau_1}k^{1-n}\|{\bf b}_2\|_\infty\|\bq- \mathring{\bq}\|_\Xi.
		\end{equation*}
	\end{corollary}
	\begin{proof}
		Let ${\bf v}=T^{-1}_{\bq}{\bf b}$ and $\mathring{\bf v}=T^{-1}_{\mathring{\bq}}{\bf b}$. Then
		\begin{equation*}
			T{\bf v}={\bf b}+\gamma \hbp,\quad T\mathring{\bf v}={\bf b}+\mathring\gamma \mathring\hbp.
		\end{equation*} 
		Taking difference, we get that
		\begin{equation*}
			T({\bf v}-\mathring{\bf v})=\mathring\gamma (\mathring\hbp-\bp)+
			(\gamma-\mathring{\gamma})\hbp.
		\end{equation*}
		As in previous proof,  we have $|\mathring\gamma|\le CR^{-1}\|{\bf b}_2\|_\infty$. On the other hand, $\mathring\hbp-\bp=M_0(\mathring\bq^{\perp}-\bq^{\perp})$. Then the results follows from Proposition \ref{prop sol T} and some simple computations.
	\end{proof}

	Now we can state the main result of this subsection.
	\begin{proposition}\label{prop reduce al}
		Under the assumption of Theorem \ref{thm u}, for sufficiently large $k$ and for each $\alpha\in\mathbb{R}$, there exists a unique $(\bq, \gamma)=(\bq(\alpha), \gamma(\alpha))$ such that
		\begin{equation}\label{reduce al}
			\left\{\begin{aligned}
				&L(\phi)+E+N(\phi)=\gamma(\alpha)\frac{\partial W}{\partial \alpha},\\
				&\int_{\R}\phi U_{j}^{p-1}Z_{j,l}\rd x=0,\ l=0,1,2, \ j=1,\dots,k.
			\end{aligned}\right.
		\end{equation}
		Moreover, the function $\phi=\phi(x;\alpha,\bq(\alpha))$ is $C^1$ in $\alpha$, and we have
		\begin{equation*}
			\|\phi\|_*\le C,\quad \|\bq\|_\Xi+R^{-1}\|\partial_\alpha\bq\|_{\Xi}\le Cd^{-\tau_2},
		\end{equation*}
		where $\tau_2\in (0,\frac{1}{2}\tau_1)$ is a small number.
	\end{proposition}
	\begin{proof}
		By previous discussion, it suffices to prove $\hbc=0$ in \eqref{eq hbc}. We only consider the case $n\ge 6$. By Lemma \ref{lem proj exp}, this is equivalent to
		\begin{equation}\label{eq red Tq}
			T\bq+{\boldsymbol{\Phi}}(\alpha,\bq)
			=\gamma \hbp,
		\end{equation}
		where the matrix $T$ is given by \eqref{def T}, and the remainder term ${\boldsymbol{\Phi}}=(\Phi_0, \Phi_1, \Phi_2)^{\T}$ satisfies
		\begin{equation*}
			\|\Phi_0\|_\infty\le Cd^{-1}k^{n-1}\log d, \ \|\Phi_0\|_{\infty'}\le C(d^{-2}k^{n-1}\log d+d^{-1-\frac{3}{2}\tau_1}k^{n-3}).
		\end{equation*}
		\begin{equation*}
			\|\Phi_1\|_\infty\le Cd^{-2\tau_1}k^{n-2},\  \|\Phi_2\|_\infty\le Cd^{-1-2\tau_1}k^{n-3}.
		\end{equation*}
		By Proposition \ref{prop sol T}, equation \eqref{eq red Tq} is equivalent to
		\begin{equation*}
			\bq=-T_{\bq}^{-1}({\boldsymbol{\Phi}}(\bq)):=\mathcal{F}(\bq).
		\end{equation*}
		
		Using estimates in Proposition \ref{prop sol T}, it is straightforward to check that $\|\FR(\bq)\|_\Xi\le Cd^{-\tau_2}$ by our choice of $\tau_2$. On the other hand, for $\bq, \mathring{\bq}\in \Xi$, by Corollary \ref{coro Tq -1}, we have
		\begin{equation*}
			\begin{aligned}
				\|\FR(\bq)-\FR(\mathring{\bq})\|_\Xi\le \|T_{\mathring\bq}^{-1}({\boldsymbol{\Phi}}(\mathring\bq)-{\boldsymbol{\Phi}}(\bq))\|_\Xi+\|(T_{\mathring\bq}^{-1}-T_{\bq}^{-1})\boldsymbol{\Phi}(\bq)\|_\Xi\le\frac{1}{2}\|\bq- \mathring{\bq}\|_\Xi.
			\end{aligned}
		\end{equation*}
		Therefore by $\FR$ is a contraction mapping on $\Xi$ and the existence then follows.
		
		The differentiability of $\bq(\alpha)$ can be easily verified by the implicit function theorem, and the estimate of $\|\partial_\alpha\bq\|_{\Xi}$ is similar to the proof of Corollary \ref{coro Tq -1}. We refer the readers to Lemma 5.5 in \cite{DWY} and omit the details here.
	\end{proof}
	
	\subsection{Variational reduction}\label{var redu}
	~~~~In this subsection, our purpose is to complete the proof of Theorem \ref{thm main}. We will apply the variational reduction to solve $\gamma(\alpha)=0$ in \eqref{reduce al}. 
	
	We introduce the reduced energy functional
	\begin{equation}\label{def J al}
		\JR(\alpha):=I(W+\phi).
	\end{equation}
	Then it is easy to see that $\JR(\alpha)$ is of class $C^1$ and $2\pi$-periodic in $\alpha$.
	
	The next lemma shows that it is enough to find critical points of $\JR(\alpha)$. The proof is nowadays standard so we omit it. We refer the readers to Lemma 6.2 in \cite{AWY} and \cite{DWY}.
	\begin{lemma}\label{lem J'=0}
		If $\alpha_0$ is a critical point of $\JR(\alpha)$, then $\gamma(\alpha_0)=0$, hence the function 
		\begin{equation*}
			u(x)=U(x;\alpha_0,\bq(\alpha_0))+\phi(x;\alpha_0,\bq(\alpha_0))
		\end{equation*}
		is a solution of \eqref{eq u}.
	\end{lemma}

	\begin{proof}[Proof of Theorem \ref{thm main}]
		Since $\JR(\alpha)$ is $2\pi$ periodic and of class $C^1$, it has at
		least two critical points (maximum and minimum points) in $[0,2\pi)$. Therefore, Theorem \ref{thm u} follows from Lemma \ref{lem J'=0}, and Theorem \ref{thm main} is proved.
	\end{proof}

	\section{Non-degeneracy Results}\label{Pf non-degen}
	~~~~~~In this section, we establish the proof of Theorem \ref{thm nondegen} and Corollary \ref{thm nondegen sym}.
	\subsection{Revisit the existence problem}
	In this subsection, we briefly revisit the existence problem for \eqref{eq u} when $K(x)$ is radially symmetric and give a slightly different proof of Theorem A, so that we can obtain extra estimates which are necessary for the proof of non-degeneracy.
	
	We place the center of bubbles at $Q_{j,0}$ for $j=1,\dots, k$, which is defined by
	\begin{equation}\label{def Qj 0}
		Q_{j,0}=(r\cos\theta_{j}, r\sin\theta_{j},\mathbf{0}),
	\end{equation}
	where $\theta_{j}=\frac{2\pi(j-1)}{k}$. We define
	\begin{equation*}\label{def UR WR}
		\UR_j(x):=U_{Q_{j,0},\Lambda_0}(x), \text{ and } \WR_0(x):=\sum_{j=1}^k\UR_j(x).
	\end{equation*}
	Here we require
	\begin{equation}\label{range r la}
		r\in (R-\frac{1}{\mu^{\bar\theta}}, R+\frac{1}{\mu^{\bar\theta}}), \quad |\Lambda_0-\Lambda|\le C_1
	\end{equation}
	for some small $\bar\theta>0$, and some fixed constant $C_1$, where we recall $R$ and $\Lambda$ are determined by \eqref{R0 La 1} and $\eqref{R0 La 2}$. Then as before, $u=\WR_0+\phi$ solves \eqref{eq u} is and only if
	\begin{equation*}
		\LR(\phi)+\ER+\NR(\phi)=0,
	\end{equation*}
	where $\LR, \ER, \NR$ are defined in the same way as \eqref{def L}, \eqref{def E} and \eqref{def N} respectively.
	
	By a slight abuse of notation, in this section, we still use $\| \ \|_*$ and $\|\ \|_{**}$ to denote norms defined by \eqref{def of n1} and \eqref{def of n2}, in which the points $Q_{j}$ is replaced by $Q_{j,0}$ for $j=1, \dots, k$. In addition, we also introduce the following norm:
	\begin{equation}\label{def of n star}
		\|h\|_{\star}:=\sup \mathcal{V}_1(x)^{-1}|h(x)|,
	\end{equation}
	where
	\begin{equation}\label{def of V1}
		\mathcal{V}_1(x)=\sum_{i=1}^{k}\frac{d^{-(n-2)}}{1+|x-Q_{i,0}|^{3}}\chi_{\{x\in \omei_i\}}+\sum_{i=1}^{k}
		\frac{d^{-\max\{p,2\}\sigma}}{1+|x-Q_{i,0}|^{n+1-\max\{p,2\}\sigma}}\chi_{\{x\in \Omega'\}}.
	\end{equation}
	
	Let $\ZR_{j,0}=\frac{\partial\UR_j}{\partial \Lambda_0}$ and $\ZR_{j,1}=\frac{\partial\UR_j}{\partial r}$. Then we have the following result, whose proof is the same as Proposition \ref{prop nonline proj} so we omit it.
	\begin{proposition}\label{prop nonline proj r}
		Let $n\geq 5$. Assume that ${K}$ satisfies $(\mathbf{K}_0)$. Then there exists a positive integer $k_0$ such that for all $k\ge k_0$, there exist a unique function $\phi\in W^{2,2}(\R)\cap H_s$ and a unique multiplier $(c, c')$ such that
		\begin{equation}\label{eq phi proj r}
			\left\{\begin{aligned}
				&\LR(\phi)+\ER+\NR(\phi)=c\sum_{j=1}^k\UR_j^{p-1}\ZR_{j,0}
				+c'\sum_{j=1}^k\UR_j^{p-1}\ZR_{j,1},\\
				&\int_{\R}\phi \UR_{j}^{p-1}\ZR_{j,l}\rd x=0, \ j=1,\dots,k, \ l=0,1.
			\end{aligned}\right.
		\end{equation}
		Moreover, there exists a large constant $C$ such that
		\begin{equation}\label{est grad phi}
			\|\phi\|_*+\|\nabla\phi\|_{\star}\le C.
		\end{equation}
	\end{proposition}
	
	With Proposition \ref{prop nonline proj r}, one can then proceed as in \cite{WY} by performing a mini-max argument to find $r$ and $\Lambda_0$ in the given range \eqref{range r la} such that $u=\WR_0+\phi$ is a genuine solution to \eqref{eq u}.
	
	However, we can say more if $K$ satisfies $(\mathbf{K}_2)$. In fact, instead of variational reduction method used in \cite{WY}, we can apply direct method as we did in the proof of Theorem \ref{thm u}. To be more precise, we multiply \eqref{eq phi proj r} by $\ZR_{j,0}$ and $\ZR_{j,1}$, and integrate in $\R$. Then by results of \eqref{E Zj 1} and \eqref{E Zj 0}, after some simplification, we obtain that $c=c'=0$ if and only if
	\begin{align}
		&c_0\mu^{-m}(\gamma_2\Lambda_0^{2-m}(r-\mu r_0)+\frac{\gamma_4\Lambda_0^{-m}}{r})+\gamma_1\Lambda_0^{2-n}(2r)^{1-n}\sum_{l=1}^{k-1} \sin^{2-n}(l\theta)=O(\frac{\log d}{d^{n+1}}+\frac{1}{d^{n+\tau}k^{2}}).\label{r La0}\\
		&\frac{1}{2}\gamma_1\Lambda_0^{2-n}(2r)^{2-n}\sum_{l=1}^{k-1}\sin^{2-n}(l\theta)-c_0\mu^{-m}\gamma_3\Lambda_0^{-m}=O(\frac{\log d}{d^n}).\label{r La1}
	\end{align}
	This explains why we introduce the balancing condition \eqref{R0 La 1} and \eqref{R0 La 2}. By a simple analysis, one can easily show that
	\begin{equation}\label{range r la new}
		|r-R|<\frac{1}{\mu d},\quad  |\Lambda_0-\Lambda|\le \frac{1}{d}.
	\end{equation}
	This improves \eqref{range r la}.
	
	\subsection{Proof of Theorem \ref{thm nondegen}} 
	~~~~~~Let $u$ be the $k$-bump solution constructed in previous subsection. Let $\LB$ be the linearized operator around $u$, that is, 
	\begin{equation}\label{def LB}
		\LB(\eta)=\Delta\eta+p\hat{K}u^{p-1}\eta.
	\end{equation}
	Then the space $\bar{\mathcal{Z}}_u$ defined by \eqref{def bar Zv} is $(2n-3)$-dimensional and is spanned by the following functions:
	\begin{equation}\label{def Z al}
		z_{\alpha}=
		\begin{cases}
			\frac{1}{r}\left( x_1\frac{\partial}{\partial x_{\alpha+2}} u(x) - x_{\alpha+2} \frac{\partial}{\partial x_1} u(x)\right),\quad & 1 \leq \alpha \leq n-2, \\
			\frac{1}{r}\left(x_2 \frac{\partial}{\partial x_{\alpha-n+4}} u(x) - x_{\alpha-n+4} \frac{\partial}{\partial x_2} u(x)\right),\quad & n-1 \leq \alpha \leq 2n-4, \\
			\frac{1}{r}\left(x_1 \frac{\partial}{\partial x_2} u(x) - x_2 \frac{\partial}{\partial x_1} u(x)\right),\quad & \alpha=2n-3.
		\end{cases}
	\end{equation}
	We will prove Theorem \ref{thm nondegen} by showing the following:
	\begin{theorem}\label{thm nondegen 1}
		Under the same assumption as Theorem \ref{thm nondegen}. Let $\LB$ be defined by \eqref{def LB}. If $\eta$ is a bounded solution to the equation 	\begin{equation*}
			\LB(\eta)=0,
		\end{equation*}
		then $\eta$ is a linear combination of the functions $z_\alpha(x)$, for $\alpha=1,\dots, 2n-3$.
	\end{theorem}
	
	Before we start the proof of Theorem \ref{thm nondegen 1}, following the procedure of \cite{AMW} and \cite{MW}, we introduce some approximate kernels of $\LB$. For writing convenience, we let 
	$	\pi_{\alpha} := \frac{\partial \phi}{\partial x_\alpha}, \ \alpha=1, 2, \dots, n.$ 
	Then by Proposition \ref{prop nonline proj r}, we have
	$\|\pi_{\alpha}\|_{\star}\le C$.
	Next, for $j=1, \dots, k$, we define 
	\begin{equation}\label{def Zj approx}
		\begin{aligned}
			Z_j^1 &= \cos\theta_j \left(\frac{\partial \UR_j}{\partial x_1}\right) + \sin\theta_j \left(\frac{\partial \UR_j}{\partial x_2}\right), \\
			Z_j^2 &= \sin\theta_j \left(\frac{\partial \UR_j}{\partial x_1}\right) - \cos\theta_j \left(\frac{\partial \UR_j}{\partial x_2}\right)+\frac{1}{kr}(x_1 \pi_2 -x_2 \pi_1), \\
			Z_{j}^{\alpha} &=\frac{\partial \UR_j}{\partial x_{\alpha}}+
			\frac{2}{rk}\bigg[\cos\theta_j(x_1 \pi_{\alpha}-x_{\alpha}\pi_1)+\sin\theta_j(x_2 \pi_{\alpha}-x_{\alpha}\pi_2)\bigg],\quad 3 \leq \alpha \leq n.
		\end{aligned}
	\end{equation}
	
	Though the definition of $Z_{j}^{\alpha}$ looks complicated, it has the advantage that the following decompositions hold:
	\begin{equation}\label{z al decom}
		z_{\alpha} = 
		\begin{cases}
			\sum\limits_{\substack{j=1}}^{k} \cos \theta_{j} Z_{j}^{(\alpha+2)}, & 1 \leq \alpha \leq n-2, \\
			\sum\limits_{\substack{j=1}}^{k} \sin \theta_{j} Z_{j}^{(\alpha-n+4)}, & n-1 \leq \alpha \leq 2n-4, \\
			\sum\limits_{\substack{j=1}}^{k} Z_{j}^{2}, & \alpha=2 n-3.
		\end{cases}
	\end{equation}
	
	Moreover, since $\phi$ and $\pi_{\alpha}$ are small enough by \eqref{est grad phi}, they play no role in all estimates we need to use, and many results in Section \ref{reduction} still hold. Now we divide the proof into following steps:
	
	\noindent\textbf{Step 1:}
	By an orthogonal decomposition, it is easy to see that we only need to show $\eta=0$ if $\eta$ satisfies the following additional condition:
	\begin{equation}\label{eta1}
		\int_{\mathbb{R}^n}{\eta}u^{p-1}z_\alpha \rd x=0\quad\hbox{for $1\leq \alpha\leq 2n-3$}.
	\end{equation} 
	We decompose ${\eta}$ using approximate kernels $Z_{j}^i$ and write it in the following form
	\begin{equation}\label{l11}
		{\eta}		=\sum_{i=0}^{n}\sum_{j=1}^{k}c_{ij}Z_{j}^{i}+\eta^{\bot},
	\end{equation}
	where $\mathbf{c}_i =(c_{i1},\dots, c_{ik})^\T$ is a $k$-dimensional vector and $\eta^{\bot}$ satisfies
	\begin{equation*}
		\int_{\mathbb{R}^n}\eta^{\bot} \UR_{j}^{p-1}Z_{j}^{i} \rd x=0\quad\text{ for $i=0, \dots, n$ and $j=1, \dots, k$},
	\end{equation*}
	Then we need to prove $\eta^{\bot}=0$ and $\mathbf{c}_i =0$ for $i=0, \dots, n$.\\
	
	\noindent\textbf{Step 2:} By \eqref{l11} and the fact $\LB({\eta})=0$, we have
	\begin{equation}\label{L eta perp}
		\LB(\eta^{\bot})=
		-\sum_{j=1}^{k}\sum_{i=0}^{n}c_{ij}\LB(Z_{j}^{i}).
	\end{equation}
	By Lemma \ref{coro norm L}, we have $\sum_{j=1}^k\|\LB(Z_j^i)\|_{**}\leq C$ for $i=0, \dots, n$.  On the other hand, the operator $\LB$ satisfies the same type of a priori estimate as $L$ in Proposition \ref{prop line theo}, so we obtain that 
	\begin{equation}\label{est eta le c}
		\|\eta^{\bot}\|_{*}\leq C\sum_{i=0}^{n}\|\mathbf{c}_i\|_{\infty}.
	\end{equation}

	\noindent\textbf{Step 3:}	Substituting \eqref{l11} into \eqref{eta1}, we find that for $1\leq {\alpha}\leq 2n-3$,
	\begin{equation}\label{eta2}
		\int_{\mathbb{R}^n}\sum_{i=0}^{n}\sum_{j=1}^{k}c_{ij}Z_{j}^{i}z_{\alpha} u^{p-1}\rd x
		=-\int_{\mathbb{R}^n}\eta^{\bot} z_{\alpha} u^{p-1}\rd x.
	\end{equation}
	We introduce three vectors in $\mathbb{R}^k$:
	\begin{equation*}
		\begin{bmatrix}
			\mathbf{1}
		\end{bmatrix}= \begin{bmatrix}
			1 \\
			1 \\
			\vdots \\
			1
		\end{bmatrix},
		\quad
		\begin{bmatrix}
			\mathbf{cos}
		\end{bmatrix} = \begin{bmatrix}
			1 \\
			\cos\theta_2 \\
			\vdots \\
			\cos\theta_k 
		\end{bmatrix},
		\quad
		\begin{bmatrix}
			\mathbf{sin}
		\end{bmatrix} = \begin{bmatrix}
			0 \\
			\sin\theta_2 \\
			\vdots \\
			\sin\theta_k 
		\end{bmatrix},
	\end{equation*}
	where  we recall
	$\theta_j =\frac{2\pi(j-1)}{k}$. Then we have following proposition:
	\begin{proposition}\label{prop c condi}
		The system \eqref{eta2} reduces to $(2n-3)$ linear conditions on vectors $\mathbf{c}_i$ for $i=0, \dots, n$. To be more specific, for $\alpha=3, \dots, n$, we have
		\begin{equation}\label{ci3}
			\begin{aligned}
				&\mathbf{c}_{\alpha} \cdot \begin{bmatrix} \mathbf{cos} \end{bmatrix} = t_{\alpha-2} + O(d^{2-n}) \mathcal{L}_{\alpha-2}(\mathbf{c}_\alpha), \\
				&\mathbf{c}_{\alpha} \cdot \begin{bmatrix} \mathbf{sin} \end{bmatrix} = t_{\alpha+n-4} + O(d^{2-n}) \mathcal{L}_{\alpha+n-4}(\mathbf{c}_\alpha),
			\end{aligned}
		\end{equation}
		and
		\begin{equation}\label{ci1}
			\mathbf{c}_2 \cdot \begin{bmatrix} \mathbf{1} \end{bmatrix} = t_{2\alpha-3} + O(d^{2-n}) \mathcal{L}_{2\alpha-3}
			\begin{bmatrix} 
				\mathbf{c}_0 \\
				\mathbf{c}_1 \\ \mathbf{c}_2 \end{bmatrix},
		\end{equation}
		where $\mathcal{L}_{l}:~\mathbb{R}^{k}\to\mathbb{R}$  for $l=1, \dots, 2n-4$, and $\mathcal{L}_{2n-3}:~\mathbb{R}^{3k}\to\mathbb{R}$ are linear functions whose coefficients are constants uniformly bounded as $k$ tends to infinity. Moreover,  for $l=1, \cdots, 2n-3$, $t_l$ are numbers that satisfy
		\begin{equation}\label{t bound}
			|t_{l}|\leq Ck d^{2-n} \|\eta^\perp\|_*.
		\end{equation}
	\end{proposition}
	\begin{proof}
		First we notice that due to the smallness of $\phi$, it is easy to check that the integrals  $\int_{\mathbb{R}^n}u^{p-1}Z_{j}^{i} Z_{m}^{l}\rd x$ satisfies the same type of estimates as \eqref{almost ortho}, so we will not list everything in details here. Now we take $\alpha=2n-3$ in \eqref{eta2}. Since $z_{2n-3}=\sum_{l=1}^{k}Z^2_l$, we have that
		\begin{equation*}
			\begin{aligned}
				-\int_{\mathbb{R}^n}\eta^{\bot} (\sum_{l=1}^{k}Z^2_l) u^{p-1}\rd x
				=&\sum_{l=1}^k\left( \sum_{j\neq l}\sum_{i=0}^2c_{ij}
				\int_{\mathbb{R}^n}Z_{j}^{i}Z^2_l u^{p-1}\rd x+
				\sum_{i=0}^2c_{il}
				\int_{\R}Z_{l}^{i}Z^2_l u^{p-1}\rd x\right)\\
				=&\sum_{l=1}^k\left( c_{2l}\tilde{\gamma}
				+(c_{0l}+c_{1l}+c_{2l})O(d^{2-n})\right),
			\end{aligned}
		\end{equation*}
		where $\tilde{\gamma}=\Lambda^2\int_{\mathbb{R}^n}U^{p-1}\left( \partial_1 U\right)^{2}\rd x$. Set 
		\begin{equation*}
			t_{2n-3}:=-\frac{1}{\tilde{\gamma}}\int_{\mathbb{R}^n}\eta^{\bot} z_{2n-3} u^{p-1} \rd x=-\frac{1}{\tilde{\gamma}}\int_{\mathbb{R}^n}\eta^{\bot}(\sum_{l=1}^kZ_l^2) u^{p-1} \rd x.
		\end{equation*}
		Then a straightforward computation gives that $|t_{2n-3}|\leq Ckd^{2-n}\|\eta^{\bot}\|_{*}$. Thus we established	\eqref{ci1}. Similarly, by setting
		\begin{equation*}
			t_{\alpha}:=-\frac{1}{\tilde{\gamma}}\int_{\mathbb{R}^n}\eta^{\bot} z_{\alpha} u^{p-1} \rd x,\quad \alpha=1,\dots,2n-4,
		\end{equation*}
		we can deduce \eqref{ci3} in the same way. The details are omitted.
	\end{proof}
	
	\noindent\textbf{Step 4:}	To estimate the size of $\mathbf{c}_i$  more accurately, we multiply both sides of \eqref{L eta perp} by $Z^i_{j}$ for $i=0,1, \dots, n$ and $j=1, \dots, k$. Then we obtain the following linear system
	\begin{equation}\label{matr}
		M \begin{bmatrix}
			\mathbf{c}_0 \\
			\vdots \\
			\mathbf{c}_n
		\end{bmatrix}
		=
		- \begin{bmatrix}
			\mathbf{r}_0 \\
			\vdots \\
			\mathbf{r}_n
		\end{bmatrix},
	\end{equation}
	where
	$
	M=\left(\int_{\mathbb{R}^n} \LB(Z_j^i) Z_m^l \rd x\right)_{ij, lm}
	$ and
	\begin{equation*}
		\mathbf{r}_{i} = 
		\begin{bmatrix}
			\int_{\mathbb{R}^n} \LB(\eta^{\perp}) Z^i_{ 1} \rd x\\
			\vdots \\
			\int_{\mathbb{R}^n} \LB(\eta^{\perp}) Z^i_{ k}\rd x
		\end{bmatrix}
		=\begin{bmatrix}
			\int_{\mathbb{R}^n} \LB(Z^i_{ 1})\eta^{\perp}  \rd x\\
			\vdots \\
			\int_{\mathbb{R}^n} \LB(Z^i_{ k})\eta^{\perp} \rd x
		\end{bmatrix}.
	\end{equation*}
	By computation similar to \eqref{est L Zjl phi} and \eqref{est L Zj0 phi}, we can obtain
	\begin{equation}\label{est ri}
		\|\mathbf{r}_{i}\|_{\infty}
		\le C\|\eta^{\perp}\|_{*}\left\{
		\begin{aligned}
			&d^{-n-3}\log d, &n\ge 7,\\
			&d^{-2(n-2)}, &n=5,6,
		\end{aligned}\right.
		i\neq 0, \quad
		\|\mathbf{r}_{0}\|_{\infty}
		\le C\|\eta^{\perp}\|_{*}\left\{
		\begin{aligned}
			&d^{-n-2}, &n\ge 7,\\
			&d^{-8}\log d, &n=6,\\
			&d^{-6}, &n=5.
		\end{aligned}\right.
	\end{equation}
	
	Concerning the solvability of the matrix equation \eqref{matr}, we establish the following proposition, whose proof is left to Appendix \ref{pf prop sol M}.
	\begin{proposition}\label{prop sol M}
		There exists an integer $k_0$  such that for all $k>k_0$, system \eqref{matr} is solvable. Moreover, the solution has the form
		\begin{equation}\label{sol c1 2}
			\begin{bmatrix}
				\mathbf{c}_0 \\
				\mathbf{c}_1 \\
				\mathbf{c}_2
			\end{bmatrix}
			=
			\begin{bmatrix}
				\mathbf{v}_0 \\
				\mathbf{v}_1 \\
				\mathbf{v}_2
			\end{bmatrix}
			+
			s_1
			\begin{bmatrix}
				\mathbf{0} \\
				\mathbf{0} \\
				\mathbf{1} 
			\end{bmatrix}
		\end{equation}
		and
		\begin{equation}\label{sol c al}
			\mathbf{c}_\alpha = \mathbf{v}_\alpha + 
			s_{\alpha 1}
			\begin{bmatrix}
				\mathbf{cos}
			\end{bmatrix}
			+ s_{\alpha 2}
			\begin{bmatrix}
				\mathbf{sin}
			\end{bmatrix}, \quad \alpha = 3, \ldots, n,
		\end{equation}
		where $s_1, s_{\alpha 1}, s_{\alpha 2} \in \mathbb{R}$ and $\mathbf{v}_{i}$ are fixed vectors that satisfy
		\begin{equation}\label{v0 1 2 infty}
			\begin{aligned}
				&\|\mathbf{v}_{0}\|_{\infty}\leq C\mu^{m}(k\|\mathbf{r}_{0}\|_\infty+\frac{ k}{\mu}\|\mathbf{r}_{1}\|_\infty+\frac{\log k}{\mu}\|\mathbf{r}_{2}\|_\infty),\\
				&\|\mathbf{v}_{1}\|_{\infty}\leq C\mu^{m}(\mu^{-1}\|\mathbf{r}_{0}\|_\infty+k\|\mathbf{r}_{1}\|_\infty+\log k \|\mathbf{r}_{2}\|_\infty),\\
				&\|\mathbf{v}_{2}\|_{\infty}\leq C\mu^{m}(\frac{\log k}{\mu}\|\mathbf{r}_{0}\|_\infty+\log k\|\mathbf{r}_{1}\|_\infty+\mu^2\|\mathbf{r}_{2}\|_\infty),
			\end{aligned}
		\end{equation}
		and
		\begin{equation}\label{v al infty}
			\|\mathbf{v}_{\alpha}\|_{\infty}\leq C\mu^{m+2}\|\mathbf{r}_\alpha\|_\infty, \quad
			\alpha=3, \dots, n.
		\end{equation}
		Consequently, in view of \eqref{est ri}, we have
		\begin{equation}\label{est v le eta}
			\|\mathbf{v}_{i}\|_{\infty}\leq	Cd^{-\tau_3}\|\eta^{\perp}\|_{*},\quad i=1,\dots,n,
		\end{equation}
		where $\tau_3>0$ is a small number that satisfies $\tau_3<\frac{3m-6}{m}$ if $n=5$ and $\tau_3<\frac{2(2m+2-2n)}{m}$ if $n\ge 6$.
	\end{proposition}
	
	\noindent\textbf{Step 5:} Now we are ready to prove $\eta=0$ by above estimates.
	We first claim that $\mathbf{c}_{i}$ is unique and satisfies
	\begin{equation}\label{ci infty up}
		\|\mathbf{c}_{i}\|_{\infty}\leq Cd^{-\tau_3}\|\eta^\perp\|_{*},\quad i=1,\dots,n.
	\end{equation}
	By  Proposition \ref{prop sol M}, we only need to prove the uniqueness of
	$s_1$, $s_{\alpha1}$, and $s_{\alpha2}$, $\alpha=3,\dots,n$.
	In fact, by \eqref{ci1} and \eqref{sol c1 2}, we get that
	\begin{equation*}
		s_1=\frac{t_{2\alpha-3}}{k}+\frac{1}{k}
		O(d^{2-n}) \mathcal{L}_{2\alpha-3}
		\left(\begin{bmatrix} 
			\mathbf{v}_0 \\ \mathbf{v}_1 \\ \mathbf{v}_2
		\end{bmatrix}
		+s_1\begin{bmatrix}
			\mathbf{0} \\
			\mathbf{0} \\
			\mathbf{1} 
		\end{bmatrix}\right).
	\end{equation*}
	Then we see that $s_1$ is uniquely determined. Moreover, by \eqref{t bound} and \eqref{est v le eta}, we have
	\begin{equation*}
		|s_1|\le Cd^{2-n}\|\eta^{\perp}\|_*+Cd^{2-n-\tau_3}\|\eta^{\perp}\|_*
		\le Cd^{2-n}\|\eta^{\perp}\|_*.
	\end{equation*} 
	This and \eqref{est v le eta} prove \eqref{ci infty up} for $i=1,2$. Similarly, we can use \eqref{ci3} and \eqref{sol c al} to derive uniqueness and the bound for $s_{\alpha1}$, and $s_{\alpha2}$, $\alpha=3,\dots,n$, hence prove \eqref{ci infty up} for $i=3,\dots,n$.
	
	Now it follows from  \eqref{est eta le c} and \eqref{ci infty up} that 
	\begin{equation*}
		\|\eta^{\bot}\|_{*}\leq C\sum_{i=1}^{n}\|\mathbf{c}_i\|_{\infty}\le Cd^{-\tau_3}\|\eta^\perp\|_{*},
	\end{equation*}
	which indicates that
	\begin{equation*}
		\mathbf{c}_{i}=\mathbf{0},\ i=1,\dots,n, \text{ and } \ \eta^\perp=0
	\end{equation*}
	provided that $k$ is large enough. As a result, ${\eta}=0$ by \eqref{l11}. Thus, we complete the proof of Theorem \ref{thm nondegen 1} and Theorem \ref{thm nondegen}.
	
	\subsection{Proof of Corollary \ref{thm nondegen sym}}
	We prove Corollary \ref{thm nondegen sym} by proving the following proposition:
	\begin{proposition}\label{prop nondegen sym 1}
		Under the same assumption as Corollary \ref{thm nondegen sym}. Let $H_s$ be defined by \eqref{def Hs} and $\LB$ be defined by \eqref{def LB}.
		If $\eta\in H_s$ is a bounded solution to the equation 	\begin{equation*}
			\LB(\eta)=0,
		\end{equation*}
		then $\eta=0$.
	\end{proposition}
	\begin{proof}
		Since $\eta\in H_s$, we can decompose
		\begin{equation}\label{eta sym decom}
			{\eta}
			=c\sum_{j=1}^{k}Z_{j}^{0}+c'\sum_{j=1}^{k}Z_{j}^{1}+\eta^{\bot},
		\end{equation}
		for some numbers $c$ and $c'$, where $\eta^{\bot}$ satisfies
		\begin{equation*}
			\int_{\mathbb{R}^n}\eta^{\bot} \UR_{j}^{p-1}Z_{j}^{i} \rd x=0\quad\text{ for $i=0, \dots, n$ and $j=1, \dots, k$}.
		\end{equation*}	
		Then we have
		\begin{equation}\label{L eta perp sym}
			\LB(\eta^{\perp})=-c\sum_{j=1}^{k}\LB(Z_{j}^{0})-c'\sum_{j=1}^{k}\LB(Z_{j}^{1}).
		\end{equation}
		as in \eqref{est eta le c}, we can deduce that
		\begin{equation}\label{est eta le c sym}
			\|\eta^{\perp}\|_*\le C(|c|+|c'|).
		\end{equation}
		
		On the other hand, by testing \eqref{L eta perp sym} against $Z_{1}^0$ and $Z_1^1$, we can obtain
		\begin{equation*}
			\begin{pmatrix}
				a_0 & b_0\\
				b_0 &  f_0\\
			\end{pmatrix}
			\begin{pmatrix}
				c\\
				c'
			\end{pmatrix}=-\begin{pmatrix}
				\int_{\R}\LB(\eta^{\perp})Z_1^0\rd x\\
				\int_{\R}\LB(\eta^{\perp})Z_1^1\rd x
			\end{pmatrix},
		\end{equation*}
		where $a_0, b_0, f_0$ are given by \eqref{a 0}, \eqref{b nu} and \eqref{f nu} respectively. We remark that the computation in Appendix \ref{pf prop sol M} is under condition $(\mathbf{K}_2)$, and the estimates under condition $(\mathbf{K}_0)$ may be not so sharp. However, we still have
		\begin{equation*}
			a_0=\frac{n-m-2}{2}\gamma_3\Lambda^{-m-2}c_0\mu^{-m}(1+o(1)), \ f_0=\mu^{-m}\zeta(1+o(1)),\ 	b_0=\mu^{-m}o(1).
		\end{equation*}
		Then by \eqref{est ri}, it is easy to see that
		\begin{equation}\label{est c le eta sym}
			|c|+|c'|\le Cd^{-\tau_4}\|\eta^{\perp}\|_*
		\end{equation}
		for some $\tau_4>0$ small enough. Combining \eqref{est eta le c sym} and \eqref{est c le eta sym}, we deduce that $c=c'=0$ and $\eta^{\perp}=0$. This proves the theorem.   
	\end{proof}
	
	\section*{Acknowledgment}
	\smalltext{ The research of J.C. Wei is supported by National Key R\&D Program of China 2022YFA1005602, and
		Hong Kong General Research Funds “New frontiers in singularity formations of nonlinear partial differential equations” and “On Fujita equation in
		critical and supercritical regime”. This paper was completed when H. Yang was visiting
		the Department of Mathematics of the Chinese University of Hong Kong. He would like to thank them for hospitality. He also thanks Central South University and School of Mathematics and Statistic for their support in his life during his visit.}
	
	\appendix
	\section{Proof of Lemma \ref{lem proj exp}}\label{pf lem proj exp}
	~~~~\normalsize We only prove the lemma for $n\ge 6$. Recall that $Q_j=(R+f_j)\vn_j+g_j\vt_j$. For $j\neq l$, we compute
	\begin{equation*}\label{Qj-Ql 2}
		\begin{aligned}
			|Q_j-Q_l|^2=&2R^2(1-\cos(\theta_j-\theta_l))
			+2R\left((f_l+f_j)(1-\cos(\theta_j-\theta_l))-(g_l-g_j)\sin(\theta_j-\theta_l)\right)\\
			&+(g_j-g_l)^2+4g_jg_l\sin^2\frac{\theta_j-\theta_l}{2}+O(\frac{d_{jl}}{R^2})\\
			=&d_{jl}^2\left(1-\frac{2R}{d_{jl}^2}(g_l-g_j)\sin(\theta_j-\theta_l)+O(\frac{d^{-\frac{3}{2}\tau_1}}{R^2})\right).
		\end{aligned}
	\end{equation*}
	As a result, for $s>0$, we have
	\begin{equation}\label{Qj-Ql -s}
		|Q_j-Q_l|^{-s}=d_{jl}^{-s}+sR d_{jl}^{-s-2}(g_l-g_j)\sin(\theta_j-\theta_l)+O(d_{jl}^{-s-\frac{3}{2}\tau_1}R^{-2}).
	\end{equation}
	
	For the normal direction, using \eqref{Qj-Ql -s} and the fact that $\bq\in\Xi$, we can expand
	\begin{equation}\label{proj Zj 1 1}
		\begin{aligned}
			&\Lambda_j^{-\frac{n-2}{2}}\sum_{l\neq j}\Lambda_l^{-\frac{n-2}{2}}
			\frac{Q_j-Q_l}{|Q_j-Q_l|^n}\cdot \vn_j\\
			=&\Lambda^{2-n}(2R)^{1-n}\sum_{l\neq j} |\sin^{2-n}((j-l)\theta)|
			-\frac{n-2}{2}\Lambda^{1-n}(2R)^{1-n}\sum_{l\neq j} |\sin^{2-n}((j-l)\theta)|(\lambda_l+\lambda_j)\\
			&-\frac{n}{2}\Lambda^{2-n}(2R)^{-n}\sum_{l\neq j} \sin(2(j-l)\theta) |\sin^{-n}((j-l)\theta)| g_j\\
			&+\frac{n-2}{2}\Lambda^{2-n}(2R)^{-n}\sum_{l\neq j} \sin(2(j-l)\theta) |\sin^{-n}((j-l)\theta)|  g_l+O(d^{-n-1}).
		\end{aligned}
	\end{equation}
	
	Since $R=\mu r_0+R_0$ with $|R_0|\le\frac{C}{\mu}$, we have
	\begin{equation}\label{proj Zj 1 2}
		\begin{aligned}
			(\gamma_2\Lambda_j^{{2}-m}(|Q_j-\mu r_0|)
			\frac{Q_j}{|Q_j|}
			+\gamma_4\Lambda_j^{-m}\frac{Q_j}{|Q_j|^2})\cdot\vn_j
			=&(\gamma_2\Lambda^{2-m}R_0+\frac{\gamma_4\Lambda^{-m}}{R})
			+\gamma_2\Lambda_j^{{2}-m}(f_j+\frac{g_j^2}{2R})\\
			&+(\gamma_2(2-m)\Lambda^{2-m}R_0-\frac{m\gamma_4\Lambda^{-m}}{R})\lambda_j+O(\mu^{-1}d^{-2}).
		\end{aligned}
	\end{equation}
	Plugging \eqref{proj Zj 1 1} and \eqref{proj Zj 1 2} into \eqref{E Zj 1}, and using the relation \eqref{R0 La 1}, we obtain
	\begin{equation}\label{proj Zj 1}
		\begin{aligned}
			\int_{\R}EZ_{j,1}\rd x
			=&c_0\mu^{-m}\gamma_2\Lambda^{{2}-m}f_j+O(d^{1-n-2\tau_1}k^{-1}).
		\end{aligned}
	\end{equation}

	For the tangential direction, using \eqref{Qj-Ql -s} again, we obtain
	\begin{equation}\label{proj Zj 2 1}
		\begin{aligned}
			&\Lambda_j^{-\frac{n-2}{2}}	\sum_{l\neq j}\Lambda_l^{-\frac{n-2}{2}}
			\frac{Q_j-Q_l}{|Q_j-Q_l|^n}\cdot \vt_j\\
			=&-\frac{n-2}{4}\Lambda^{1-n}(2R)^{1-n}\sum_{l\neq j}\sin(2(j-l)\theta) |\sin^{2-n}((j-l)\theta)|(\lambda_l+\lambda_j)\\
			&+\Lambda^{2-n}(2R)^{-n}\sum_{l\neq j} |\sin^{-n}((j-l)\theta)| 
			\left( n-1-(n-2)\sin^2((j-l)\theta)\right) (g_l-g_j)\\
			&+2\Lambda^{2-n}(2R)^{-n}\sum_{l\neq j} |\sin^{2-n}((j-l)\theta)|g_j+O(d^{-n-1}k^{-2}).
		\end{aligned}
	\end{equation}
	Similar to \eqref{proj Zj 1 2}, we also have
	\begin{equation}\label{proj Zj 2 2}
		\begin{aligned}
			(\gamma_2\Lambda_j^{{2}-m}|Q_j-\mu r_0|
			\frac{Q_j}{|Q_j|}
			+\gamma_4\Lambda_j^{-m}\frac{Q_j}{|Q_j|^2})\cdot\vt_j=(\gamma_2\Lambda^{2-m}R_0+\frac{\gamma_4\Lambda^{-m}}{R})\frac{g_j}{R}+O(\mu^{-2}d^{-2\tau_1}).
		\end{aligned}
	\end{equation}
	
	Plugging \eqref{proj Zj 2 1} and \eqref{proj Zj 2 2} into \eqref{E Zj 2}, and using the relation \eqref{R0 La 1} again, we obtain
	\begin{equation}\label{proj Zj 2}
		\begin{aligned}
			\int_{\R}EZ_{j,2}\rd x	=
			&\gamma_1\Lambda^{2-n}(2R)^{-n}\sum_{l\neq j} |\sin^{-n}((j-l)\theta)| 
			\left( n-1-(n-2)\sin^2((j-l)\theta)\right) (g_l-g_j)\\
			&-\frac{n-2}{4}\gamma_1\Lambda^{1-n}(2R)^{1-n}\sum_{l\neq j}\sin(2(j-l)\theta) |\sin^{2-n}((j-l)\theta)|\lambda_l+O(d^{-n-2\tau_1}k^{-2}).
		\end{aligned}
	\end{equation}
	
	Finally, for the projection along $Z_{j,0}$, using \eqref{Qj-Ql -s}, we compute
	\begin{equation*}
		\begin{aligned}
			&\Lambda_j^{-\frac{n-2}{2}}
			\sum_{l\neq j}\Lambda_l^{-\frac{n-2}{2}}
			\frac{1}{|Q_j-Q_l|^{n-2}}\\
			=&\Lambda^{2-n}(1-\frac{n-2}{2}\Lambda^{-1}(\lambda_j+\lambda_l))
			\sum_{l\neq j}\left(d_{jl}^{2-n}+(n-2)Rd_{jl}^{-n}(g_l-g_j)\sin(\theta_j-\theta_l)\right)+O(d^{-n-\frac{3}{2}\tau_1}k^{-2})\\
			=&\Lambda^{2-n}(2R)^{2-n}\sum_{l\neq j}|\sin^{2-n}(j-l)\theta|
			-\frac{n-2}{2}\Lambda^{1-n}(2R)^{2-n}\sum_{l\neq j}|\sin^{2-n}(j-l)\theta|(\lambda_j+\lambda_l)\\
			&+\frac{n-2}{2}\Lambda^{2-n} (2R)^{1-n}\sum_{l\neq j}|\sin^{-n}(j-l)\theta|\sin(\theta_j-\theta_l)(g_l-g_j)+O(d^{-n-\frac{3}{2}\tau_1}k^{-2}),
		\end{aligned}
	\end{equation*}
	and
	\begin{equation*}
		c_0\mu^{-m}\gamma_3\Lambda_j^{-m}=c_0\mu^{-m}\gamma_3\Lambda^{-m}-mc_0\mu^{-m}\gamma_3\Lambda^{1-m}\lambda_j+O(d^{-n-2\tau_1}k^{-2}).
	\end{equation*}
	Therefore by \eqref{E Zj 0} and the relation \eqref{R0 La 2}, we obtain
	\begin{equation}\label{proj Zj 0}
		\begin{aligned}
			\int_{\R}EZ_{j,0}\rd x=&-\frac{n-2}{2}\gamma_1\Lambda^{1-n}(2R)^{2-n}\sum_{l\neq j}|\sin^{2-n}(j-l)\theta|\lambda_l\\
			&+(m-\frac{n-2}{2})\gamma_1\Lambda^{1-n}(2R)^{2-n}\sum_{l\neq j}|\sin^{2-n}(j-l)\theta|\lambda_j+\hat{I}\\
			&+\frac{n-2}{2}\gamma_1\Lambda^{1-n} (2R)^{1-n}\sum_{l\neq j}|\sin^{-n}(j-l)\theta|\sin(2(j-l)\theta)g_l+O(d^{-n-\frac{3}{2}\tau_1}k^{-2}+d^{-n-1}{\log d}).
		\end{aligned}
	\end{equation}
	
	By Lemma \ref{lem L+N Zj 12} and Lemma \ref{lem L+N Zj 0}, $\int_{\R} (L(\phi)+N(\phi))Z_{j,l}\rd x$ are small enough for all $j=1,\dots, k$ and $l=0,1,2$, so  \eqref{proj exp} follows from \eqref{proj Zj 1}, \eqref{proj Zj 2} and \eqref{proj Zj 0}.

	\section{Proof of  Lemma \ref{lem sol T'}}\label{pf lem sol T'}
	~~~~We first introduce some properties of circulant matrices. It is well-known that the eigenvalues of a $(k\times k)$-circulant matrix \( Y=\text{Cir}[
	y_0, y_1, \ldots, y_{k-1}] \) of dimension \( k \times k \) are given explicitly by the formula
	\begin{equation}\label{eigencir}
		\eta_\nu = \sum_{l=0}^{k-1} y_l e^{\frac{2\pi \nu}{k} i l}, \quad \nu = 0, \ldots, k-1, 
	\end{equation}
	and their corresponding normalized eigenvectors are given by
	\begin{equation*}
		E_\nu = k^{-\frac{1}{2}} 
		(1,
		e^{\frac{2\pi \nu}{k} i},
		\cdots
		e^{\frac{2\pi \nu}{k} i (k-1)})^\T, \quad \nu = 0, \ldots, k-1.
	\end{equation*}
	Moreover, $Y$ can be diagonalized as 
	$Y=PD_YP^{\top}$,
	where $D_Y$ is the diagonal matrix
	$D_Y=\text{diag}(\eta_0, \eta_1, \ldots, \eta_{k-1})$
	and $P$ is the $k \times k$ invertible matrix defined by
	\begin{equation}\label{def P}
		P = \left[ \begin{array}{c|c|c|c} E_0 & E_1 & \ldots & E_{k-1} \end{array} \right].
	\end{equation}
	For more background knowledge on circulant matrices, we refer the readers to \cite{KS}.
	
	Now we start to analysis the matrix $T$. Let $\lambda_{1,\nu}, \lambda_{2,\nu}, \lambda_{3,\nu}, \nu=0,1,\dots, k-1$, be the eigenvalues of $A_1, A_2$ and $A_3$ respectively. Then we have
	\begin{equation}\label{la 1 nu}
		\lambda_{1,\nu}=\sum_{l=1}^{k-1} \frac{\cos(2\nu l\theta)}{\sin^{n-2}(l\theta)}
		=2(\frac{k}{\pi})^{n-2}g''(\frac{2\pi\nu}{k})+O(k^{n-4}\log k),
	\end{equation}
	\begin{equation}\label{la 2 nu}
		\lambda_{2,\nu}=i\sum_{l=1}^{k-1} \frac{\sin(2\nu l\theta)\sin(2l\theta)}{\sin^{n}(l\theta)}
		=i4(\frac{k}{\pi})^{n-1}g'(\frac{2\pi\nu}{k})+\nu O(k^{n-4}\log k).
	\end{equation}
	\begin{equation}\label{la 3 nu}
		\lambda_{3,\nu}=\sum_{l=1}^{k-1} \frac{\cos(2\nu l\theta)-1}{\sin^{n}(l\theta)}
		-\frac{n-2}{n-1}\sum_{l=1}^{k-1} \frac{\cos(2\nu l\theta)-1}{\sin^{n-2}(l\theta)}
		=-2(\frac{k}{\pi})^{n}g(\frac{2\pi\nu}{k})
		+\nu^2O(k^{n-4}).
	\end{equation}
	In particular,
	\begin{equation*}
		\lambda_{2,0}=\lambda_{3,0}=0.
	\end{equation*}
	
	Let $\mathcal{P}=\text{diag}(P, P)$,
	where $P$ is defined in \eqref{def P}. Then we can write
	\begin{equation*}
		\mathcal{P}^\T T'\mathcal{P}
		=D':=\begin{pmatrix} 
			c_1D_1+c_2I&  D_2 \\
			\frac{1}{2} D_2 &  c_3D_3 \\
		\end{pmatrix},
	\end{equation*}
	where $D_1=\text{diag}(\lambda_{1,0},\dots,\lambda_{1,k-1})$, and $D_2, D_3$ are defined in the same way.
	
	Let ${\bf y}=\mathcal{P}^\T {\bf v}$ and ${\bf r}=\mathcal{P}^\T {\bf b}$. Then \eqref{eq T'v} is equivalent to
	\begin{equation}\label{eq D'y}
		D'{\bf y}={\bf r}+\gamma \mathcal{P}^\T\bp_0',\quad y_{2,0}=0,
	\end{equation}
	which can be decomposed into the following $k$ equations:
	\begin{equation}\label{eq nu=0}
		\begin{pmatrix} 
			c_1\lambda_{1,0}+c_2&  0 \\
			0&  0\\
		\end{pmatrix}
		\begin{pmatrix}
			y_{0,0}\\
			0\\
		\end{pmatrix}
		=\begin{pmatrix}
			r_{0,0}\\
			r_{2,0}\\
		\end{pmatrix}
		+\gamma\begin{pmatrix}
			0\\
			k^{\frac{1}{2}}\\
		\end{pmatrix},
	\end{equation}
	\begin{equation}\label{eq nu neq 0}
		\begin{pmatrix} 
			c_1\lambda_{1,\nu}+c_2&  \lambda_{2,\nu} \\
			-\frac{1}{2} \lambda_{2,\nu} &  c_3\lambda_{3,\nu} \\
		\end{pmatrix}
		\begin{pmatrix}
			y_{1,\nu}\\
			y_{2,\nu}\\
		\end{pmatrix}
		=\begin{pmatrix}
			r_{1,\nu}\\
			r_{2,\nu}\\
		\end{pmatrix},\quad \nu=1,\dots,k-1.
	\end{equation}
	
	For \eqref{eq nu=0} to be solvable, we need $\gamma=k^{-\frac{1}{2}}r_{2,0}$, and the solution is
	\begin{equation}\label{y 0 0}
		y_{0,0}=\frac{r_{0,0}}{c_1\lambda_{1,0}+c_2}.
	\end{equation}
	
	To solve \eqref{eq nu neq 0}, we denote $\mathcal{D}_\nu=\text{det}(D_\nu)=(c_1\lambda_{1,\nu}+c_2)c_3\lambda_{3,\nu} -\frac{1}{2}\lambda_{2,\nu}^2$.
	Then by \eqref{la 1 nu}, \eqref{la 2 nu} and \eqref{la 3 nu}, we can compute
	\begin{equation}\label{D nu}
		\begin{aligned}
			\mathcal{D}_\nu
			=&8(\frac{k}{\pi})^{2n-2}\left(
			(g'(\frac{2\pi\nu}{k}))^2-\frac{n-1}{n-2}
			g(\frac{2\pi\nu}{k})
			g''(\frac{2\pi\nu}{k})\right)
			\\
			&+\frac{8(2m+2-n)(n-1)}{(n-2)^2}(\frac{k}{\pi})^{2n-2}g(0)g(\frac{2\pi\nu}{k}).
		\end{aligned}
	\end{equation}
	Since $2m+2-n>0$, by assumption \eqref{condi g}, we find that there exist positive constants $C_1, C_2$, which are independent of $k$ and $\nu$, such that
	\begin{equation}\label{est det nu}
		C_1\nu^2 k^{2n-4}\le	\mathcal{D}_\nu \le C_2\nu^2 k^{2n-4},\ \nu=1,\dots,k.
	\end{equation}
	Now the solution of \eqref{eq nu neq 0} is given by
	\begin{equation}\label{y nu}
		y_{0,\nu}=\frac{(c_3\lambda_{3,\nu})r_{0,\nu}
			-\lambda_{2,\nu}r_{2,\nu}}{\mathcal{D}_\nu},\quad
		y_{2,\nu}=\frac{\frac{1}{2}\lambda_{2,\nu}r_{0,\nu}
			+(c_1\lambda_{1,\nu}+c_2)r_{2,\nu}}{\mathcal{D}_\nu}.
	\end{equation}
	
	Using \eqref{y 0 0} and \eqref{y nu}, and transforming $({\bf y}, {\bf r})$ back to $({\bf v}, {\bf b})$, we can obtain
	\begin{equation*}
		\begin{aligned}
			v_{0,j}
			=&k^{-\frac{1}{2}}\frac{r_{0,0}}{c_1\lambda_{1,0}+c_2}
			+k^{-\frac{1}{2}}\left(\sum_{\nu=1}^{k-1}\frac{(c_3\lambda_{3,\nu})r_{0,\nu}
				-\lambda_{2,\nu}r_{2,\nu}}{\mathcal{D}_\nu}e^{\frac{2\pi i}{k}j\nu}\right)\\
			=&k^{-1}\sum_{l=0}^{k-1}
			\left((\frac{1}{c_1\lambda_{1,0}+c_2}
			+\sum_{\nu=1}^{k-1}\frac{c_3\lambda_{3,\nu}}{\mathcal{D}_\nu}e^{\frac{2\pi i}{k}(l+j)\nu})b_{0,l}
			+(\sum_{\nu=1}^{k-1}\frac{-\lambda_{2,\nu}}{\mathcal{D}_\nu}	e^{\frac{2\pi i}{k}(l+j)\nu})b_{2,l}
			\right).
		\end{aligned}
	\end{equation*}
	Therefore
	\begin{equation*}
		\begin{aligned}
			|v_{0,j}|\le& \frac{C}{|c_1\lambda_{1,0}+c_2|}\|{\bf b}_0\|_{\infty}
			+C\sum_{\nu=1}^{k-1}
			\Big|\frac{c_3\lambda_{3,\nu}}{\mathcal{D}_\nu}\Big|\|{\bf b}_0\|_{\infty'}
			+C\sum_{\nu=1}^{k-1}
			\Big|\frac{\lambda_{2,\nu}}{\mathcal{D}_\nu}\Big|
			\|{\bf b}_2\|_{\infty}\\
			\le &\frac{C}{R}\sum_{l=2}^{k}\frac{1}{\sin^{n-2}((l-1)\theta)}\|{\bf b}_0\|_{\infty}
			+\frac{C}{R}\sum_{\nu=1}^{k-1}
			\frac{\nu^2 k^{n-2}}{\nu^2 k^{2n-4}}\|{\bf b}_0\|_{\infty'}
			+C\sum_{\nu=1}^{k-1}
			\frac{\nu k^{n-2}}{\nu^2 k^{2n-4}}
			\|{\bf b}_2\|_{\infty}\\
			\le& Cd^{-1}k^{1-n}\|{\bf b}_0\|_{\infty}+Cd^{-1}k^{2-n}\|{\bf b}_0\|_{\infty'}+Ck^{2-n}\log k\|{\bf b}_2\|_{\infty}.      
		\end{aligned}
	\end{equation*}
	
	Similarly, we have
	\begin{equation*}
		\begin{aligned}
			|v_{2,j}|=&k^{-1}\Big|\sum_{l=0}^{k-1}
			\left((\sum_{\nu=1}^{k-1}\frac{\lambda_{2,\nu}}{2\mathcal{D}_\nu}e^{\frac{2\pi i}{k}(l+j)\nu})b_{0,l}	+(\sum_{\nu=1}^{k-1}\frac{(c_1\lambda_{1,\nu}+c_2)}{\mathcal{D}_\nu}	e^{\frac{2\pi i}{k}(l+j)\nu})b_{2,l}			\right)\Big|\\
			\le & Ck^{2-n}\log k\|{\bf b}_0\|_{\infty}+Cdk^{3-n}\|{\bf b}_2\|_{\infty}.     
		\end{aligned}
	\end{equation*}

	Moreover, for $j\neq t$, we have
	\begin{equation*}
		\begin{aligned}
			v_{2,j}-v_{2,t}=&k^{-1}\sum_{l=0}^{k-1}
			\Bigg(\sum_{\nu=1}^{k-1}\frac{\lambda_{2,\nu}}{2\mathcal{D}_\nu}e^{\frac{2\pi i}{k}(l+j)\nu}(1-e^{\frac{2\pi i}{k}(t-j)\nu})b_{0,l}\\	&+\sum_{\nu=1}^{k-1}\frac{(c_1\lambda_{1,\nu}+c_2)}{\mathcal{D}_\nu}	e^{\frac{2\pi i}{k}(l+j)\nu}(1-e^{\frac{2\pi i}{k}(t-j)\nu})b_{2,l}
			\Bigg).
		\end{aligned}
	\end{equation*}
	Without loss of generality, we can assume $|t-j|\le \frac{k}{2}$. Then we have
	\begin{equation*}
		\begin{aligned}
			|v_{2,j}-v_{2,t}|
			\le &Ck^{2-n}\sum_{\nu=1}^{k-1}\frac{\sin^2((t-j)\nu\theta)}{\nu}\|{\bf b}_0\|_{\infty'}
			+CRk^{2-n}\sum_{\nu=1}^{k-1}\frac{\sin^2((t-j)\nu\theta)}{\nu^2}\|{\bf b}_2\|_{\infty}\\
			\le & C\log(1+|t-j|)k^{2-n}\|{\bf b}_0\|_{\infty'}+Cdk^{2-n}|t-j|\|{\bf b}_2\|_{\infty}.     
		\end{aligned}
	\end{equation*}
	This proves \eqref{est v2-}, and the proof of the lemma is completed.

	\section{Proof of  Proposition \ref{prop sol M}}\label{pf prop sol M}
	~~~~We first give a comprehensive analysis of the matrix $M$.
	By symmetry, if $l\in\{3,\cdots,n\}$ and $m\in\{0,1,2\}$, then $\int_{\mathbb{R}^n} \LB(Z_{i}^l) Z^m_{j} \, \rd x = 0$ for any $i, j = 1, \ldots, k$, so the matrix $M$ has the form
	\begin{equation*}
		M = \begin{bmatrix} M_1 & 0 \\ 0 & M_2 \end{bmatrix}
	\end{equation*}
	where $M_1$ is a $(3k\times 3k)$-dimensional matrix and $M_2$ is an $((n-2)k\times (n-2)k)$-dimensional matrix.
	To be more specific, $M_1$ has following form:
	\begin{equation*}
		M_1 = \begin{bmatrix} 
			A & B & C\\
			B^{\text{T}} & F & E\\
			C^{\text{T}} & E^{\text{T}} & G.
		\end{bmatrix}
	\end{equation*}
	where
	\begin{equation*}
		\begin{aligned}
			A =& \left(\int_{\mathbb{R}^n} \LB(Z_i^0) Z_j^0 \rd x\right)_{i, j=1, \dots, k},\quad B = \left(\int_{\mathbb{R}^n} \LB(Z_i^0) Z_j^1 \rd x\right)_{i, j=1, \dots, k},\\
			C =& \left(\int_{\mathbb{R}^n} \LB(Z_i^0) Z_j^2 \rd x\right)_{i, j=1, \dots, k},\quad F = \left(\int_{\mathbb{R}^n} \LB(Z_i^1) Z_j^1 \rd x\right)_{i, j=1, \dots, k},\\
			E =& \left(\int_{\mathbb{R}^n} \LB(Z_i^1) Z_j^2 \rd x\right)_{i, j=1, \dots, k},\quad G = \left(\int_{\mathbb{R}^n} \LB(Z_i^2) Z_j^2 \rd x\right)_{i, j=1, \dots, k}
		\end{aligned}
	\end{equation*}
	are all $(k\times k)$-dimensional matrices. Similarly, since
	\[
	\int_{\mathbb{R}^n} \LB(Z^l_{i}) Z^m_{ j} \rd x = 0, \quad \text{if} \quad l \neq m, \quad l,m = 3, \ldots, n,
	\]
	the matrix \( M_2 \) has the form
	$M_2 = \text{diag}(H_3,H_4,\dots, H_n)$,
	where
	\begin{equation*}
		H_\alpha=\left(\int_{\mathbb{R}^n} \LB(Z_i^{\alpha}) Z_j^{\alpha} \rd x\right)_{i, j=1, \dots, k}
	\end{equation*}
	are also $(k\times k)$ matrices, $\alpha = 3, \ldots, n$. Moreover, following similar arguments as presented in \cite{AMW, MW}, we can conclude that $A, B, C, E, F, G$ and $H_{\alpha}$ are all circulant matrices. 
	
	Next we describe the entries of these matrices, together with their eigenvalues. The computations of these entries is actually a paraphrase of computations in Lemma \ref{lem E Zj 12} and Lemma \ref{lem E Zj 0}, and we only need to apply \eqref{E3}, \eqref{E4} and \eqref{E5} repeatedly, so omit the details.\\
	
	\noindent\text{\bf The matrix $A$}.
	The matrix  $A$ defined by $A_{ij}=\int_{\mathbb{R}^n} \LB(Z_i^0) Z_j^0 \rd x$ is symmetric. We have
	\begin{equation}\label{A11}
		A_{11}=\frac{n-m-2}{2}\gamma_3\Lambda^{-m-2}c_0\mu^{-m}+O(d^{-n}),
	\end{equation}
	and for $j>1$,
	\begin{equation}\label{A1j}
		A_{1j}=\frac{n-2}{4}\gamma_1
		\Lambda^{-n}(2r)^{2-n}\sin^{2-n}\frac{\theta_j}{2}+O(d_j^{2-n}d^{-2}).
	\end{equation}
	Then the eigenvalues of $A$ are
	\begin{equation}\label{a nu}
		\begin{aligned}
			a_\nu
			=&\frac{n-2m-2}{2}\gamma_3\Lambda^{-m-2}c_0\mu^{-m}
			+\frac{n-2}{2}\gamma_1
			\Lambda^{-n}(2r)^{2-n}
			(\frac{k}{\pi})^{n-2}g''(\frac{2\pi\nu}{k})+O(d^{2-n}k^{-2}\log k).
		\end{aligned}
	\end{equation}
	In particular, by \eqref{R0 La 2}, we have
	\begin{equation}\label{a 0}
		a_0=\frac{n-m-2}{2}\gamma_3\Lambda^{-m-2}c_0\mu^{-m}+O(d^{2-n}k^{-2}\log k).
	\end{equation}
	
	\noindent\text{\bf The matrix $B$}.
	The matrix  $B$ defined by $B_{ij}=\int_{\mathbb{R}^n} \LB(Z_i^0) Z_j^1 \rd x$ is symmetric. We have
	\begin{equation}\label{B11}
		\begin{aligned}
			B_{11}
			=&-(m-2)c_0\mu^{-m}(|Q_1|-\mu r_0)\gamma_2\Lambda^{1-m}-\frac{mc_0\mu^{-m}}{|Q_j|}\gamma_4\Lambda^{-1-m}\\
			&+\frac{n-2}{2}\gamma_1\Lambda^{1-n}(2r)^{1-n}
			\sum_{j\neq 1 }\sin^{n-2}\frac{\theta_j}{2}+O(\mu^{-m-2})+O(d^{-n-1}), 
		\end{aligned}
	\end{equation}
	and for $j>1$,
	\begin{equation}\label{B1j}
		B_{1j}=-\frac{n-2}{4}\gamma_1\Lambda^{1-n}(2r)^{1-n}\sin^{2-n}\frac{\theta_j}{2}+O(d_j^{1-n}d^{-3}).
	\end{equation}
	Then the eigenvalues of $B$ are
	\begin{equation}\label{b nu}
		\begin{aligned}
			b_\nu=&B_{11}-\frac{n-2}{4}\gamma_1\Lambda^{1-n}(2r)^{1-n}
			(\frac{k}{\pi})^{n-2}g''(\frac{2\pi\nu}{k})+O(\mu^{-m-2})+O(d^{-n-1})+O(d^{1-n}k^{-3}).
		\end{aligned}
	\end{equation}
	They satisfies
	\begin{equation}\label{est b nu}
		|b_\nu|\le C\mu^{-m-1}.
	\end{equation}
	
	\noindent\text{\bf The matrix $C$}.
	The matrix $C$ defined by $C_{ij}=\int_{\mathbb{R}^n} \LB(Z_i^0) Z_j^2 \rd x$ is anti-symmetric. We have
	\begin{equation}\label{C 1j}
		C_{1j}=-\frac{n-2}{4}\gamma_1\Lambda^{1-n}(2r)^{1-n}\sin{\theta_j}\sin^{-n}\frac{\theta_j}{2}+O(d_j^{1-n}d^{-3}), \ j=2,\dots ,k.
	\end{equation}
	By our choice of $Z_j^2$, we can use \eqref{z al decom} to deduce that $C_{11}=-\sum_{j=2}^{k}C_{1j}=0$. This implies that $C$ has a zero eigenvalue $c_0=\sum^{k}_{j=1} C_{1j}=0$. For $\nu= 1, \dots, k-1$, the other $k-1$ eigenvalues $c_{\nu}$ of matrix $C$ can be written as
	\begin{equation}
		\begin{aligned}\label{c nu}
			c_\nu			=i\sum_{j=2}^kC_{1j}\sin\nu\theta_j
			=-(n-2)i\gamma_1\Lambda^{1-n}(2r)^{1-n}(\frac{k}{\pi})^{n-1}g'(\frac{2\pi\nu}{k})
			+O(\nu d^{1-n}k^{-3}).
		\end{aligned}
	\end{equation}
	Since $|g'(x)|\le Cx$, we have
	\begin{equation}\label{est c nu}
		|c_\nu|\le C\mu^{-m-1}\nu.
	\end{equation}
	
	\noindent\text{\bf The matrix $F$}.
	The matrix $F$ defined by $ F_{ij}= \int_{\mathbb{R}^n} \LB(Z_i^1) Z_j^1 \rd x$ is symmetric.  We have
	\begin{equation}\label{F11}
		F_{11}=\mu^{-m}\zeta+O(\mu^{-m-2-\tau}),
	\end{equation}
	where $\zeta:=\frac{c_{0}m(m-1)}{p+1}
	\int_{\mathbb{R}^n}|y|^{m-2}U^{p+1}(y)\rd y.$ For $j\neq 1$, we have
	\begin{equation}\label{F1j}
		F_{1j}=\gamma_1\Lambda^{2-n} (2 r)^{-n}(\sin^{-n}\frac{\theta_{j}}{2}+(n-2)\sin^{2-n}\frac{\theta_{j}}{2})+O(d^{-3}d_j^{1-n}).
	\end{equation}
	Then it is easy to see that for $\nu=0, 1, \dots, k-1,$ the eigenvalues $f_{\nu}$ of matrix $F$ satisfies
	\begin{equation}\label{f nu}
		\begin{aligned}
			f_{\nu} 
			=&\mu^{-m}\zeta+O(\mu^{-m-2-\tau}).
		\end{aligned}
	\end{equation}
	
	\noindent\text{\bf The matrix $E$}.
	The matrix $E$ defined by $E_{ij}=\int_{\mathbb{R}^n} \LB(Z_i^1) Z_j^2 \rd x$ is anti-symmetric. Similar to the matrix $C$, we have
	\begin{equation}\label{E1j}
		\begin{aligned}
			E_{1j}=\frac{n-2}{2}\gamma_1\Lambda^{2-n}(2r)^{-n}\sin\theta_j\sin^{-n}\frac{\theta_j}{2}+O(d_j^{1-n}d^{-3}), \quad j\neq 1,
		\end{aligned}
	\end{equation}
	and 	$E_{11}=\sum_{j=2}^kE_{1j}=0$. Then $E$ has a zero eigenvalue $e_0=0$. For $\nu=1, \dots, k-1,$ the other $k-1$ eigenvalues $e_\nu$ of $E$ are given by
	\begin{equation}\label{e nu}
		\begin{aligned}
			e_{\nu}=&i2({n-2})\gamma_1\Lambda^{2-n}(2r)^{-n}(\frac{k}{\pi})^{n-1}g'(\frac{2\pi\nu}{k})+\nu O(\mu^{-n}k^{n-4}\log k)
			+O(d^{-2-n}g'(\frac{2\pi\nu}{k})).
		\end{aligned}
	\end{equation}
	Then we obtain
	\begin{equation}\label{est e nu}
		|e_\nu|\le C \mu^{-m-2}\nu.
	\end{equation}
	
	\noindent\text{\bf The matrix $G$}.
	The matrix $G$ defined by $G_{ij}=\int_{\mathbb{R}^n} \LB(Z_i^2) Z_j^2 \rd x$ is symmetric. We have
	\begin{equation}\label{G1j}
		G_{1j}=	\gamma_1\Lambda^{2-n} (2r)^{-n}\left((1-n)\sin^{-n}\frac{\theta_j}{2}+(n-2)\sin^{2-n}\frac{\theta_j}{2}\right)+O(d^{-3}d_j^{1-n}),\ j=2,\dots,k.
	\end{equation}
	Moreover, by our choice of $Z_j^2$ again, we have $G_{11}=-\sum_{j=2}^{k}G_{1j}$, which implies $g_0=0$.
	For $\nu= 1, \dots, k-1$, the other $k-1$ eigenvalues $g_{\nu}$ of matrix $G$ can be written as
	\begin{equation}\label{g nu}
		\begin{split}
			g_\nu&=\gamma_1\Lambda^{2-n} (2r)^{-n}(n-1)\left(2(\frac{k}{\pi})^{n}g(\frac{2\pi\nu}{k})+O(\nu^2k^{n-4}\log k)\right).
		\end{split}
	\end{equation}
	Then we have
	\begin{equation}\label{est g nu}
		c\nu^2\mu^{-m-2}\le	g_\nu \le C\nu^2\mu^{-m-2},\ \nu=1,\dots,k-1.
	\end{equation}
	for some small constant $c$ and large constant $C$.
	
	\noindent\text{\bf The matrix $H_\alpha$}.
	The  matrix $H_\alpha$ defined by	$	H_{\alpha,ij}=\int_{\mathbb{R}^n} \LB(Z_i^{\alpha}) Z_j^{\alpha} \rd x$ is symmetric. For $j\neq 1$, similar to the computation of $A$ and $G$, we have that
	\begin{equation}\label{H al 1j}
		H_{\alpha,1j}=2\gamma_1\Lambda^{2-n} (2r)^{-n}\sin^{-n}\frac{\theta_j}{2}+O(d_j^{1-n}d^{-3}).
	\end{equation}
	
	Moreover, for $\alpha\ge 3$, by our choice of $Z_j^\alpha$, we can use \eqref{z al decom} to deduce that
	\begin{equation}\label{sum H=0}
		\sum_{j=1}^k \cos \theta_j H_{\alpha,1j}=	\sum_{j=1}^k \sin \theta_j H_{\alpha,1j}=0.
	\end{equation}
	
	Concerning the eigenvalues of $H_\alpha$, we first notice that \eqref{sum H=0} implies that
	\begin{equation}\label{h al 1 k-1}
		h_{\alpha,1}=h_{\alpha,k-1}=\sum_{l=j}^{k}H_{\alpha,1j} \cos \theta_j =0.
	\end{equation}
	For $\nu=0,2,3, \dots, k-2$, the other $k-2$ eigenvalues $h_{\alpha,\nu}$ of $H_\alpha$ are given by
	\begin{equation}\label{h al nu}
		\begin{aligned}
			h_{\alpha,\nu}
			&=2\gamma_1\Lambda^{2-n} (2r)^{-n} 2(\frac{k}{\pi})^n\left(g(\frac{2\pi}{k})-g(\frac{2\pi\nu}{k})\right)
			+O(\mu^{-n}k^{n-4}(\nu^2+1)).
		\end{aligned}
	\end{equation}
	Then it is easy to see that there is a small constant $c$ such that
	\begin{equation}\label{est h al}
		|h_{\alpha,\nu}|\ge c(\nu^2+1)\mu^{-m-2}, \ \text{ for all } \nu=0,2,3,\dots,k-2.
	\end{equation}
	
	Now we turn to the solvability of \eqref{matr}. We only prove the part about $M_1$, since the part about $M_2$ is similar and in fact easier. Let $\mathcal{P}_1=\text{diag}(P, P, P)$,
	where $P$ is defined by \eqref{def P}. Then $M_1$ can be written as
	$$
	M_1 =  \mathcal{P}_1 {D}_{M_1} \mathcal{P}_1^{{\text{T}}},
	$$
	where
	\begin{equation}\label{def D M1}
		{D}_{M_1}=\begin{bmatrix}
			D_A & D_B & D_C\\
			D_B & D_F & D_E\\
			-D_C & -D_E & D_G
		\end{bmatrix}.
	\end{equation}
	Here $D_X$ denotes the $(k\times k)$-dimensional diagonal matrix  whose entries are given by the eigenvalues of $X$. For instance $D_A=\text{diag}(a_0,\dots, a_{k-1})$.
	
	For $\nu=0, 1, \dots, k-1$, we set
	\begin{equation*}
		\hat{D}_\nu=\begin{pmatrix}
			a_\nu & b_\nu & c_\nu \\ 
			b_\nu & f_\nu & e_\nu \\
			-c_\nu & -e_\nu & g_\nu 
		\end{pmatrix}.   
	\end{equation*}
	and $\hat{\mathcal{D}}_\nu=\text{det}(\hat{D}_\nu)$. Then $\hat{\DR}_0=0$, and if $\nu\neq 0$, similar to the computation of $\DR_\nu$ in \eqref{D nu}, we can find  positive constants $C_1, C_2$, which are independent of $k$ and $\nu$, such that
	\begin{equation}\label{est det hat nu}
		C_1\mu^{-3m-2}\nu^2\le  -\hat{\DR}_\nu \le C_2\mu^{-3m-2}\nu^2.
	\end{equation}
	
	Let
	\begin{equation*}
		\hat{D}_\nu^{-1}=\begin{pmatrix}
			\tilde a_\nu &\tilde b_\nu &\tilde c_\nu \\ 
			\tilde b_\nu &\tilde f_\nu &\tilde e_\nu \\
			-\tilde c_\nu & -\tilde e_\nu &\tilde g_\nu 
		\end{pmatrix}.   
	\end{equation*}
	Then using \eqref{est det hat nu} and estimates of eigenvalues \eqref{a nu}, \eqref{est b nu}, \eqref{est c nu}, \eqref{f nu}, \eqref{est e nu}, \eqref{est g nu}, we obtain
	\begin{equation}\label{est D nu -1}
		\begin{aligned}
			&|\tilde a_\nu|\le C\mu^m,\  |\tilde b_\nu|\le C\mu^{m-1}, \ |\tilde c_\nu|\le C\frac{\mu^{m-1}}{\nu},\\
			&|\tilde f_\nu|\le C\mu^m, \ |\tilde e_\nu|\le C\frac{\mu^{m}}{\nu}, \ |\tilde g_\nu|\le C\frac{\mu^{m+2}}{\nu^2}.
		\end{aligned}
	\end{equation}
	
	Using the same argument as the proof of  Lemma \ref{lem sol T'}, we can write $\mathbf{c}_l$ explicitly in terms of $\mathbf{r}_l$ for $l=0,1,2$. Then estimates \eqref{v0 1 2 infty} follows from \eqref{est D nu -1}.

\end{document}